\documentclass[12pt,reqno]{amsart}

\catcode`\@=11
\def\@secnumfont{\bfseries}
\def\section{\@startsection{section}{1}%
  \z@{.7\linespacing\@plus\linespacing}{.5\linespacing}%
  {\normalfont\bfseries}}
\def\subsection{\@startsection{subsection}{2}%
  \z@{.5\linespacing\@plus.7\linespacing}{-.5em}%
  {\normalfont\scshape}}
\catcode`\@=12

\numberwithin{equation}{section}

\newtheorem*{Theorem}{Theorem}
\newtheorem*{Conjecture}{Conjecture}

\catcode`\@=11
\def\sideset#1#2#3{%
  \@mathmeasure\z@\displaystyle{#3}%
  \global\setbox\@ne\vbox to\ht\z@{}\dp\@ne\dp\z@
  \setbox\tw@\box\@ne
  \@mathmeasure4\displaystyle{\copy\tw@#1}%
  \@mathmeasure6\displaystyle{#3{#2}}%
  \dimen@-\wd6 \advance\dimen@\wd4 \advance\dimen@\wd\z@
  \hbox to\dimen@{}\mathop{\kern-\dimen@\box4\box6}%
}
\catcode`\@=12

{\obeylines
\gdef\MATH{\begingroup\parindent0pt\parskip0pt plus 0.25pt\obeylines%
        \def^^M{\vskip5pt}%
        \obeyspaces\tt\small}%
}
\def\goodbreakpoint{\par\penalty-5000%
         \vrule height10pt depth2pt width0pt\leavevmode}
\def\endMATH{\endgroup}
\def\MATHphi{\leavevmode
        \hbox to 0pt{\hbox to 5.24995pt{\hss$\phi$\hss}\hss}}
\def\MATHGamma{\leavevmode
        \hbox to 0pt{\hbox to 5.24995pt{\hss$\Gamma$\hss}\hss}}
\def\MATHpi{\leavevmode
        \hbox to 0pt{\hbox to 5.24995pt{\hss$\pi$\hss}\hss}}
\def\MATHinfty{\leavevmode
        \hbox to 0pt{\hbox to 5.24995pt{\hss$\infty$\hss}\hss}}
\def\MATHhStrich{\leavevmode
        \hbox to 0pt{\hbox to 5.24995pt{\vrule height4.5pt depth-3.5pt width5.24995pt}\hss}}
\def\MATHluEck{\leavevmode
        \hbox to 0pt{\hbox to 5.24995pt{\hskip2.12497pt
         \vrule height4.5pt depth1pt width1pt
         \vrule height4.5pt depth-3.5pt width2.12498pt}\hss}}
\def\MATHruEck{\leavevmode
        \hbox to 0pt{\hbox to 5.24995pt{%
         \vrule height4.5pt depth-3.5pt width2.12497pt
         \vrule height4.5pt depth1pt width1pt
         \hskip2.12498pt}\hss}}
\def\MATHloEck{\leavevmode
        \hbox to 0pt{\hbox to 5.24995pt{\hskip2.12497pt
         \vrule height9pt depth-3.5pt width1pt
         \vrule height4.5pt depth-3.5pt width2.12498pt}\hss}}
\def\MATHroEck{\leavevmode
        \hbox to 0pt{\hbox to 5.24995pt{%
         \vrule height4.5pt depth-3.5pt width2.12497pt
         \vrule height9pt depth-3.5pt width1pt
         \hskip2.12498pt}\hss}}
\def\MATHvStrich{\leavevmode
        \hbox to 0pt{\hbox to 5.24995pt{\hskip2.12497pt
         \vtop to 0pt{\hsize1pt\vss%
                \vrule height17pt depth6pt width1pt\vskip8pt\vss\par}%
         \hskip2.12498pt}\hss}}
\def\MATHtStueck{\leavevmode
        \hbox to 0pt{\hbox to 5.24995pt{%
         \vrule height4.5pt depth-3.5pt width2.12497pt
         \vrule height4.5pt depth2pt width1pt
         \vrule height4.5pt depth-3.5pt width2.12498pt}\hss}}
\def\MATHbackslash{\leavevmode
        \hbox to 0pt{\hbox to 5.24995pt{\hss$\backslash$\hss}\hss}}
\def\MATHlbrace{\leavevmode
        \hbox to 0pt{\hbox to 5.24995pt{\hss$\{$\hss}\hss}}
\def\MATHrbrace{\leavevmode
        \hbox to 0pt{\hbox to 5.24995pt{\hss$\}$\hss}\hss}}
\def\MATHkleiner{\leavevmode
        \hbox to 0pt{\hbox to 5.24995pt{\hss$\langle$\hss}\hss}}
\def\MATHgroesser{\leavevmode
        \hbox to 0pt{\hbox to 5.24995pt{\hss$\rangle$\hss}\hss}}
\def\MATHhoch{\leavevmode
        \hbox to 0pt{\hbox to 5.24995pt{\hss$^\land$\hss}\hss}}
\def\MATHtief{\leavevmode
        \hbox to 0pt{\hbox to 5.24995pt{\hss\vrule height0pt depth.8pt width3pt\hss}\hss}}

\def\setRevDate $#1 #2 #3${#2}
\def\TeXdrawId{\setRevDate $Date: 1995/12/19 16:40:42 $ TeXdraw V2R0}
\chardef\catamp=\the\catcode`\@
\catcode`\@=11
\long
\def\centertexdraw #1{\hbox to \hsize{\hss
\btexdraw #1\etexdraw
\hss}}
\def\btexdraw {\x@pix=0             \y@pix=0
\x@segoffpix=\x@pix  \y@segoffpix=\y@pix
\t@exdrawdef
\setbox\t@xdbox=\vbox\bgroup\offinterlineskip
\global\d@bs=0
\global\t@extonlytrue
\global\p@osinitfalse
\s@avemove \x@pix \y@pix
\m@pendingfalse
\global\p@osinitfalse
\p@athfalse
\the\everytexdraw}
\def\etexdraw {\ift@extonly \else
\t@drclose
\fi
\egroup
\ifdim \wd\t@xdbox>0pt
\t@xderror {TeXdraw box non-zero size,
possible extraneous text}%
\fi
\vbox {\offinterlineskip
\pixtobp \xminpix \l@lxbp  \pixtobp \yminpix \l@lybp
\pixtobp \xmaxpix \u@rxbp  \pixtobp \ymaxpix \u@rybp
\hbox{\t@xdinclude
[{\l@lxbp},{\l@lybp}][{\u@rxbp},{\u@rybp}]{\p@sfile}}%
\pixtodim \xminpix \t@xpos  \pixtodim \yminpix \t@ypos
\kern \t@ypos
\hbox {\kern -\t@xpos
\box\t@xdbox
\kern \t@xpos}%
\kern -\t@ypos\relax}}
\def\drawdim #1 {\def\d@dim{#1\relax}}
\def\setunitscale #1 {\edef\u@nitsc{#1}%
\realmult \u@nitsc \s@egsc \d@sc}
\def\relunitscale #1 {\realmult {#1}\u@nitsc \u@nitsc
\realmult \u@nitsc \s@egsc \d@sc}
\def\setsegscale #1 {\edef\s@egsc {#1}%
\realmult \u@nitsc \s@egsc \d@sc}
\def\relsegscale #1 {\realmult {#1}\s@egsc \s@egsc
\realmult \u@nitsc \s@egsc \d@sc}
\def\bsegment {\ifp@ath
\f@lushbs
\f@lushmove
\fi
\begingroup
\x@segoffpix=\x@pix
\y@segoffpix=\y@pix
\setsegscale 1
\global\advance \d@bs by 1\relax}
\def\esegment {\endgroup
\ifnum \d@bs=0
\writetx {es}%
\else
\global\advance \d@bs by -1
\fi}
\def\savecurrpos (#1 #2){\getsympos (#1 #2)\a@rgx\a@rgy
\s@etcsn \a@rgx {\the\x@pix}%
\s@etcsn \a@rgy {\the\y@pix}}
\def\savepos (#1 #2)(#3 #4){\getpos (#1 #2)\a@rgx\a@rgy
\coordtopix \a@rgx \t@pixa
\advance \t@pixa by \x@segoffpix
\coordtopix \a@rgy \t@pixb
\advance \t@pixb by \y@segoffpix
\getsympos (#3 #4)\a@rgx\a@rgy
\s@etcsn \a@rgx {\the\t@pixa}%
\s@etcsn \a@rgy {\the\t@pixb}}
\def\linewd #1 {\coordtopix {#1}\t@pixa
\f@lushbs
\writetx {\the\t@pixa\space sl}}
\def\setgray #1 {\f@lushbs
\writetx {#1 sg}}
\def\lpatt (#1){\listtopix (#1)\p@ixlist
\f@lushbs
\writetx {[\p@ixlist] sd}}
\def\lvec (#1 #2){\getpos (#1 #2)\a@rgx\a@rgy
\s@etpospix \a@rgx \a@rgy
\writeps {\the\x@pix\space \the\y@pix\space lv}}
\def\rlvec (#1 #2){\getpos (#1 #2)\a@rgx\a@rgy
\r@elpospix \a@rgx \a@rgy
\writeps {\the\x@pix\space \the\y@pix\space lv}}
\def\move (#1 #2){\getpos (#1 #2)\a@rgx\a@rgy
\s@etpospix \a@rgx \a@rgy
\s@avemove \x@pix \y@pix}
\def\rmove (#1 #2){\getpos (#1 #2)\a@rgx\a@rgy
\r@elpospix \a@rgx \a@rgy
\s@avemove \x@pix \y@pix}
\def\lcir r:#1 {\coordtopix {#1}\t@pixa
\writeps {\the\t@pixa\space cr}%
\r@elupd \t@pixa \t@pixa
\r@elupd {-\t@pixa}{-\t@pixa}}
\def\fcir f:#1 r:#2 {\coordtopix {#2}\t@pixa
\writeps {\the\t@pixa\space #1 fc}%
\r@elupd \t@pixa \t@pixa
\r@elupd {-\t@pixa}{-\t@pixa}}
\def\lellip rx:#1 ry:#2 {\coordtopix {#1}\t@pixa
\coordtopix {#2}\t@pixb
\writeps {\the\t@pixa\space \the\t@pixb\space el}%
\r@elupd \t@pixa \t@pixb
\r@elupd {-\t@pixa}{-\t@pixb}}
\def\fellip f:#1 rx:#2 ry:#3 {\coordtopix {#2}\t@pixa
\coordtopix {#3}\t@pixb
\writeps {\the\t@pixa\space \the\t@pixb\space #1 fe}%
\r@elupd \t@pixa \t@pixb
\r@elupd {-\t@pixa}{-\t@pixb}}
\def\larc r:#1 sd:#2 ed:#3 {\coordtopix {#1}\t@pixa
\writeps {\the\t@pixa\space #2 #3 ar}}
\def\ifill f:#1 {\writeps {#1 fl}}
\def\lfill f:#1 {\writeps {#1 fp}}
\def\htext #1{\def\testit {#1}%
\ifx \testit\l@paren
\let\next=\h@move
\else
\let\next=\h@text
\fi
\next {#1}}
\def\rtext td:#1 #2{\def\testit {#2}%
\ifx \testit\l@paren
\let\next=\r@move
\else
\let\next=\r@text
\fi
\next td:#1 {#2}}

\def\textref h:#1 v:#2 {\ifx #1R%
\edef\l@stuff {\hss}\edef\r@stuff {}%
\else
\ifx #1C%
\edef\l@stuff {\hss}\edef\r@stuff {\hss}%
\else
\edef\l@stuff {}\edef\r@stuff {\hss}%
\fi
\fi
\ifx #2T%
\edef\t@stuff {}\edef\b@stuff {\vss}%
\else
\ifx #2C%
\edef\t@stuff {\vss}\edef\b@stuff {\vss}%
\else
\edef\t@stuff {\vss}\edef\b@stuff {}%
\fi
\fi}
\def\avec (#1 #2){\getpos (#1 #2)\a@rgx\a@rgy
\s@etpospix \a@rgx \a@rgy
\writeps {\the\x@pix\space \the\y@pix\space (\a@type)
\the\a@lenpix\space \the\a@widpix\space av}}
\def\ravec (#1 #2){\getpos (#1 #2)\a@rgx\a@rgy
\r@elpospix \a@rgx \a@rgy
\writeps {\the\x@pix\space \the\y@pix\space (\a@type)
\the\a@lenpix\space \the\a@widpix\space av}}
\def\arrowheadsize l:#1 w:#2 {\coordtopix{#1}\a@lenpix
\coordtopix{#2}\a@widpix}
\def\arrowheadtype t:#1 {\edef\a@type{#1}}
\def\clvec (#1 #2)(#3 #4)(#5 #6)%
{\getpos (#1 #2)\a@rgx\a@rgy
\coordtopix \a@rgx\t@pixa
\advance \t@pixa by \x@segoffpix
\coordtopix \a@rgy\t@pixb
\advance \t@pixb by \y@segoffpix
\getpos (#3 #4)\a@rgx\a@rgy
\coordtopix \a@rgx\t@pixc
\advance \t@pixc by \x@segoffpix
\coordtopix \a@rgy\t@pixd
\advance \t@pixd by \y@segoffpix
\getpos (#5 #6)\a@rgx\a@rgy
\s@etpospix \a@rgx \a@rgy
\writeps {\the\t@pixa\space \the\t@pixb\space
\the\t@pixc\space \the\t@pixd\space
\the\x@pix\space \the\y@pix\space cv}}
\def\drawbb {\bsegment
\drawdim bp
\linewd 0.24
\setunitscale {\p@sfactor}
\writeps {\the\xminpix\space \the\yminpix\space mv}%
\writeps {\the\xminpix\space \the\ymaxpix\space lv}%
\writeps {\the\xmaxpix\space \the\ymaxpix\space lv}%
\writeps {\the\xmaxpix\space \the\yminpix\space lv}%
\writeps {\the\xminpix\space \the\yminpix\space lv}%
\esegment}
\def\getpos (#1 #2)#3#4{\g@etargxy #1 #2 {} \\#3#4%
\c@heckast #3%
\ifa@st
\g@etsympix #3\t@pixa
\advance \t@pixa by -\x@segoffpix
\pixtocoord \t@pixa #3%
\fi
\c@heckast #4%
\ifa@st
\g@etsympix #4\t@pixa
\advance \t@pixa by -\y@segoffpix
\pixtocoord \t@pixa #4%
\fi}
\def\getsympos (#1 #2)#3#4{\g@etargxy #1 #2 {} \\#3#4%
\c@heckast #3%
\ifa@st \else
\t@xderror {TeXdraw: invalid symbolic coordinate}%
\fi
\c@heckast #4%
\ifa@st \else
\t@xderror {TeXdraw: invalid symbolic coordinate}%
\fi}
\def\listtopix (#1)#2{\def #2{}%
\edef\l@ist {#1 }%
\m@oretrue
\loop
\expandafter\g@etitem \l@ist \\\a@rgx\l@ist
\a@pppix \a@rgx #2%
\ifx \l@ist\empty
\m@orefalse
\fi
\ifm@ore
\repeat}
\def\realmult #1#2#3{\dimen0=#1pt
\dimen2=#2\dimen0
\edef #3{\expandafter\c@lean\the\dimen2}}
\def\intdiv #1#2#3{\t@counta=#1
\t@countb=#2
\ifnum \t@countb<0
\t@counta=-\t@counta
\t@countb=-\t@countb
\fi
\t@countd=1
\ifnum \t@counta<0
\t@counta=-\t@counta
\t@countd=-1
\fi
\t@countc=\t@counta  \divide \t@countc by \t@countb
\t@counte=\t@countc  \multiply \t@counte by \t@countb
\advance \t@counta by -\t@counte
\t@counte=-1
\loop
\advance \t@counte by 1
\ifnum \t@counte<16
\multiply \t@countc by 2
\multiply \t@counta by 2
\ifnum \t@counta<\t@countb \else
\advance \t@countc by 1
\advance \t@counta by -\t@countb
\fi
\repeat
\divide \t@countb by 2
\ifnum \t@counta<\t@countb
\advance \t@countc by 1
\fi
\ifnum \t@countd<0
\t@countc=-\t@countc
\fi
\dimen0=\t@countc sp
\edef #3{\expandafter\c@lean\the\dimen0}}
\def\coordtopix #1#2{\dimen0=#1\d@dim
\dimen2=\d@sc\dimen0
\t@counta=\dimen2
\t@countb=\s@ppix
\divide \t@countb by 2
\ifnum \t@counta<0
\advance \t@counta by -\t@countb
\else
\advance \t@counta by \t@countb
\fi
\divide \t@counta by \s@ppix
#2=\t@counta}
\def\pixtocoord #1#2{\t@counta=#1%
\multiply \t@counta by \s@ppix
\dimen0=\d@sc\d@dim
\t@countb=\dimen0
\intdiv \t@counta \t@countb #2}
\def\pixtodim #1#2{\t@countb=#1%
\multiply \t@countb by \s@ppix
#2=\t@countb sp\relax}
\def\pixtobp #1#2{\dimen0=\p@sfactor pt
\t@counta=\dimen0
\multiply \t@counta by #1%
\ifnum \t@counta < 0
\advance \t@counta by -32768
\else
\advance \t@counta by 32768
\fi
\divide \t@counta by 65536
\edef #2{\the\t@counta}}
\newcount\t@counta    \newcount\t@countb
\newcount\t@countc    \newcount\t@countd
\newcount\t@counte
\newcount\t@pixa      \newcount\t@pixb
\newcount\t@pixc      \newcount\t@pixd
\newdimen\t@xpos      \newdimen\t@ypos
\newcount\xminpix      \newcount\xmaxpix
\newcount\yminpix      \newcount\ymaxpix
\newcount\a@lenpix     \newcount\a@widpix
\newcount\x@pix        \newcount\y@pix
\newcount\x@segoffpix  \newcount\y@segoffpix
\newcount\x@savepix    \newcount\y@savepix
\newcount\s@ppix
\newcount\d@bs
\newcount\t@xdnum
\global\t@xdnum=0
\newbox\t@xdbox
\newwrite\drawfile
\newif\ifm@pending
\newif\ifp@ath
\newif\ifa@st
\newif\ifm@ore
\newif \ift@extonly
\newif\ifp@osinit
\newtoks\everytexdraw
\def\l@paren{(}
\def\a@st{*}
\catcode`\%=12
\def\p@b {
\catcode`\%=14
\catcode`\{=12  \catcode`\}=12  \catcode`\u=1 \catcode`\v=2
\def\l@br u{v  \def\r@br u}v
\catcode `\{=1  \catcode`\}=2   \catcode`\u=11 \catcode`\v=11
{\catcode`\p=12 \catcode`\t=12
\gdef\c@lean #1pt{#1}}
\def\sppix#1/#2 {\dimen0=1#2 \s@ppix=\dimen0
\t@counta=#1%
\divide \t@counta by 2
\advance \s@ppix by \t@counta
\divide \s@ppix by #1%
\t@counta=\s@ppix
\multiply \t@counta by 65536
\advance \t@counta by 32891
\divide \t@counta by 65782
\dimen0=\t@counta sp
\edef\p@sfactor {\expandafter\c@lean\the\dimen0}}
\def\g@etargxy #1 #2 #3 #4\\#5#6{\def #5{#1}%
\ifx #5\empty
\g@etargxy #2 #3 #4 \\#5#6
\else
\def #6{#2}%
\def\next {#3}%
\ifx \next\empty \else
\t@xderror {TeXdraw: invalid coordinate}%
\fi
\fi}
\def\c@heckast #1{\expandafter
\c@heckastll #1\\}
\def\c@heckastll #1#2\\{\def\testit {#1}%
\ifx \testit\a@st
\a@sttrue
\else
\a@stfalse
\fi}
\def\g@etsympix #1#2{\expandafter
\ifx \csname #1\endcsname \relax
\t@xderror {TeXdraw: undefined symbolic coordinate}%
\fi
#2=\csname #1\endcsname}
\def\s@etcsn #1#2{\expandafter
\xdef\csname#1\endcsname {#2}}
\def\g@etitem #1 #2\\#3#4{\edef #4{#2}\edef #3{#1}}
\def\a@pppix #1#2{\edef\next {#1}%
\ifx \next\empty \else
\coordtopix {#1}\t@pixa
\ifx #2\empty
\edef #2{\the\t@pixa}%
\else
\edef #2{#2 \the\t@pixa}%
\fi
\fi}
\def\s@etpospix #1#2{\coordtopix {#1}\x@pix
\advance \x@pix by \x@segoffpix
\coordtopix {#2}\y@pix
\advance \y@pix by \y@segoffpix
\u@pdateminmax \x@pix \y@pix}
\def\r@elpospix #1#2{\coordtopix {#1}\t@pixa
\advance \x@pix by \t@pixa
\coordtopix {#2}\t@pixa
\advance \y@pix by \t@pixa
\u@pdateminmax \x@pix \y@pix}
\def\r@elupd #1#2{\t@counta=\x@pix
\advance\t@counta by #1%
\t@countb=\y@pix
\advance\t@countb by #2%
\u@pdateminmax \t@counta \t@countb}
\def\u@pdateminmax #1#2{\ifnum #1>\xmaxpix
\global\xmaxpix=#1%
\fi
\ifnum #1<\xminpix
\global\xminpix=#1%
\fi
\ifnum #2>\ymaxpix
\global\ymaxpix=#2%
\fi
\ifnum #2<\yminpix
\global\yminpix=#2%
\fi}
\def\s@avemove #1#2{\x@savepix=#1\y@savepix=#2%
\m@pendingtrue
\ifp@osinit \else
\global\p@osinittrue
\global\xminpix=\x@savepix \global\yminpix=\y@savepix
\global\xmaxpix=\x@savepix \global\ymaxpix=\y@savepix
\fi}
\def\f@lushmove {\global\p@osinittrue
\ifm@pending
\writetx {\the\x@savepix\space \the\y@savepix\space mv}%
\m@pendingfalse
\p@athfalse
\fi}
\def\f@lushbs {\loop
\ifnum \d@bs>0
\writetx {bs}%
\global\advance \d@bs by -1
\repeat}
\def\h@move #1#2 #3)#4{\move (#2 #3)%
\h@text {#4}}
\def\h@text #1{\pixtodim \x@pix \t@xpos
\pixtodim \y@pix \t@ypos
\vbox to 0pt{\normalbaselines
\t@stuff
\kern -\t@ypos
\hbox to 0pt{\l@stuff
\kern \t@xpos
\hbox {#1}%
\kern -\t@xpos
\r@stuff}%
\kern \t@ypos
\b@stuff\relax}}
\def\r@move td:#1 #2#3 #4)#5{\move (#3 #4)%
\r@text td:#1 {#5}}
\def\r@text td:#1 #2{\vbox to 0pt{\pixtodim \x@pix \t@xpos
\pixtodim \y@pix \t@ypos
\kern -\t@ypos
\hbox to 0pt{\kern \t@xpos
\rottxt {#1}{\z@sb {#2}}%
\hss}%
\vss}}
\def\z@sb #1{\vbox to 0pt{\normalbaselines
\t@stuff
\hbox to 0pt{\l@stuff \hbox {#1}\r@stuff}%
\b@stuff}}
\ifx \rotatebox\@undefined
\def\rottxt #1#2{\bgroup
#2%
\egroup}
\else
\let\rottxt=\rotatebox
\fi
\ifx \t@xderror\@undefined
\let\t@xderror=\errmessage
\fi
\def\t@exdrawdef {\sppix 300/in
\drawdim in
\edef\u@nitsc {1}%
\setsegscale 1
\arrowheadsize l:0.16 w:0.08
\arrowheadtype t:T
\textref h:L v:B }
\ifx \includegraphics\@undefined
\def\t@xdinclude [#1,#2][#3,#4]#5{%
\begingroup
\message {<#5>}%
\leavevmode
\t@counta=-#1%
\t@countb=-#2%
\setbox0=\hbox{%
\includegraphics{#5}}%
\t@ypos=#4 bp%
\advance \t@ypos by -#2 bp%
\t@xpos=#3 bp%
\advance \t@xpos by -#1 bp%
\dp0=0pt \ht0=\t@ypos  \wd0=\t@xpos
\box0%
\endgroup}
\else
\let\t@xdinclude=\includegraphics
\fi
\def\writeps #1{\f@lushbs
\f@lushmove
\p@athtrue
\writetx {#1}}
\def\writetx #1{\ift@extonly
\global\t@extonlyfalse
\t@xdpsfn \p@sfile
\t@dropen \p@sfile
\fi
\w@rps {#1}}
\def\w@rps #1{\immediate\write\drawfile {#1}}
\def\t@xdpsfn #1{%
\global\advance \t@xdnum by 1
\ifnum \t@xdnum<10
\xdef #1{\jobname.ps\the\t@xdnum}
\else
\xdef #1{\jobname.p\the\t@xdnum}
\fi
}
\def\t@dropen #1{%
\immediate\openout\drawfile=#1%
\w@rps {\p@b PS-Adobe-3.0 EPSF-3.0}%
\w@rps {\p@p BoundingBox: (atend)}%
\w@rps {\p@p Title: TeXdraw drawing: #1}%
\w@rps {\p@p Pages: 1}%
\w@rps {\p@p Creator: \TeXdrawId}%
\w@rps {\p@p CreationDate: \the\year/\the\month/\the\day}%
\w@rps {50 dict begin}%
\w@rps {/mv {stroke moveto} def}%
\w@rps {/lv {lineto} def}%
\w@rps {/st {currentpoint stroke moveto} def}%
\w@rps {/sl {st setlinewidth} def}%
\w@rps {/sd {st 0 setdash} def}%
\w@rps {/sg {st setgray} def}%
\w@rps {/bs {gsave} def /es {stroke grestore} def}%
\w@rps {/fl \l@br gsave setgray fill grestore}%
\w@rps    { currentpoint newpath moveto\r@br\space def}%
\w@rps {/fp {gsave setgray fill grestore st} def}%
\w@rps {/cv {curveto} def}%
\w@rps {/cr \l@br gsave currentpoint newpath 3 -1 roll 0 360 arc}%
\w@rps    { stroke grestore\r@br\space def}%
\w@rps {/fc \l@br gsave setgray currentpoint newpath}%
\w@rps    { 3 -1 roll 0 360 arc fill grestore\r@br\space def}%
\w@rps {/ar {gsave currentpoint newpath 5 2 roll arc stroke grestore} def}%
\w@rps {/el \l@br gsave /svm matrix currentmatrix def}%
\w@rps    { currentpoint translate scale newpath 0 0 1 0 360 arc}%
\w@rps    { svm setmatrix stroke grestore\r@br\space def}%
\w@rps {/fe \l@br gsave setgray currentpoint translate scale newpath}%
\w@rps    { 0 0 1 0 360 arc fill grestore\r@br\space def}%
\w@rps {/av \l@br /hhwid exch 2 div def /hlen exch def}%
\w@rps    { /ah exch def /tipy exch def /tipx exch def}%
\w@rps    { currentpoint /taily exch def /tailx exch def}%
\w@rps    { /dx tipx tailx sub def /dy tipy taily sub def}%
\w@rps    { /alen dx dx mul dy dy mul add sqrt def}%
\w@rps    { /blen alen hlen sub def}%
\w@rps    { gsave tailx taily translate dy dx atan rotate}%
\w@rps    { (V) ah ne {blen 0 gt {blen 0 lineto} if} {alen 0 lineto} ifelse}%
\w@rps    { stroke blen hhwid neg moveto alen 0 lineto blen hhwid lineto}%
\w@rps    { (T) ah eq {closepath} if}%
\w@rps    { (W) ah eq {gsave 1 setgray fill grestore closepath} if}%
\w@rps    { (F) ah eq {fill} {stroke} ifelse}%
\w@rps    { grestore tipx tipy moveto\r@br\space def}%
\w@rps {\p@sfactor\space \p@sfactor\space scale}%
\w@rps {1 setlinecap 1 setlinejoin}%
\w@rps {3 setlinewidth [] 0 setdash}%
\w@rps {0 0 moveto}%
}
\def\t@drclose {%
\bgroup
\w@rps {stroke end showpage}%
\w@rps {\p@p Trailer:}%
\pixtobp \xminpix \l@lxbp  \pixtobp \yminpix \l@lybp
\pixtobp \xmaxpix \u@rxbp  \pixtobp \ymaxpix \u@rybp
\w@rps {\p@p BoundingBox: \l@lxbp\space \l@lybp\space
\u@rxbp\space \u@rybp}%
\w@rps {\p@p EOF}%
\egroup
\immediate\closeout\drawfile
}
\catcode`\@=\catamp

\def\DreieckBreit{\rlvec(0 1)
                           \rlvec(0.8660254037844 -0.5)
                           \rlvec(-0.8660254037844 -0.5)}
\def\DreieckSpitz{\rlvec(0.8660254037844 0.5)
                           \rlvec(0 -1)
                           \rlvec(-0.8660254037844 0.5)}

\def\ldreieck{\bsegment
  \rlvec(0.866025403784439 .5) \rlvec(0 -1)
  \rlvec(-0.866025403784439 .5)
  \savepos(0.866025403784439 -.5)(*ex *ey)
        \esegment
  \move(*ex *ey)
        }
\def\rdreieck{\bsegment
  \rlvec(0.866025403784439 -.5) \rlvec(-0.866025403784439 -.5)  \rlvec(0 1)
  \savepos(0 -1)(*ex *ey)
        \esegment
  \move(*ex *ey)
        }
\def\rhombus{\bsegment
  \rlvec(0.866025403784439 .5) \rlvec(0.866025403784439 -.5)
  \rlvec(-0.866025403784439 -.5)  \rlvec(0 1)
  \rmove(0 -1)  \rlvec(-0.866025403784439 .5)
  \savepos(0.866025403784439 -.5)(*ex *ey)
        \esegment
  \move(*ex *ey)
        }
\def\RhombusA{\bsegment
  \rlvec(0.866025403784439 .5) \rlvec(0.866025403784439 -.5)
  \rlvec(-0.866025403784439 -.5) \rlvec(-0.866025403784439 .5)
  \savepos(0.866025403784439 -.5)(*ex *ey)
        \esegment
  \move(*ex *ey)
        }
\def\RhombusB{\bsegment
  \rlvec(0.866025403784439 .5) \rlvec(0 -1)
  \rlvec(-0.866025403784439 -.5) \rlvec(0 1)
  \savepos(0 -1)(*ex *ey)
        \esegment
  \move(*ex *ey)
        }
\def\RhombusC{\bsegment
  \rlvec(0.866025403784439 -.5) \rlvec(0 -1)
  \rlvec(-0.866025403784439 .5) \rlvec(0 1)
  \savepos(0.866025403784439 -.5)(*ex *ey)
        \esegment
  \move(*ex *ey)
        }
\def\RhombusAsh{\bsegment
  \rlvec(0.866025403784439 .5) \rlvec(0.866025403784439 -.5)
  \rlvec(-0.866025403784439 -.5) \rlvec(-0.866025403784439 .5)
  \lfill f:.8
  \savepos(0.866025403784439 -.5)(*ex *ey)
        \esegment
  \move(*ex *ey)
        }
\def\RhombusBsh{\bsegment
  \rlvec(0.866025403784439 .5) \rlvec(0 -1)
  \rlvec(-0.866025403784439 -.5) \rlvec(0 1)
  \lfill f:0.2
  \savepos(0 -1)(*ex *ey)
        \esegment
  \move(*ex *ey)
        }
\def\RhombusCsh{\bsegment
  \rlvec(0.866025403784439 -.5) \rlvec(0 -1)
  \rlvec(-0.866025403784439 .5) \rlvec(0 1)
  \savepos(0.866025403784439 -.5)(*ex *ey)
        \esegment
  \move(*ex *ey)
        }
\def\StrichV{\bsegment
  \lpatt(.05 .13)
  \rlvec(.6062177827 .35)
  \move (0 0)
  \rlvec(-.6062177827 .35)
  \savepos(0 0)(*ex *ey)
        \esegment
  \move(*ex *ey)
        }
\def\StrichA{\bsegment
  \lpatt(.05 .13)
  \rlvec(.6062177827 -.35)
  \move (0 0)
  \rlvec(-.6062177827 -.35)
  \savepos(0 0)(*ex *ey)
        \esegment
  \move(*ex *ey)
        }
\def\StrichX{\bsegment
  \lpatt(.05 .13)
  \rlvec(-.6062177827 .35)
  \move (0 0)
  \rlvec(-.6062177827 -.35)
  \savepos(0 0)(*ex *ey)
        \esegment
  \move(*ex *ey)
        }
\def\StrichY{\bsegment
  \lpatt(.05 .13)
  \rlvec(.6062177827 .35)
  \move (0 0)
  \rlvec(.6062177827 -.35)
  \savepos(0 0)(*ex *ey)
        \esegment
  \move(*ex *ey)
        }

\begin{document}

\newbox\Adr
\setbox\Adr\vbox{
\centerline{\Small\uppercase{C.~Krattenthaler$^\dagger$}}
\vskip12pt
\centerline{Institut f\"ur Mathematik der Universit\"at Wien,}
\centerline{Strudlhofgasse 4, A-1090 Wien, Austria.}
\centerline{E-mail: {\tt\footnotesize kratt@ap.univie.ac.at}}
\centerline{WWW: \footnotesize\tt http://www.mat.univie.ac.at/People/kratt}
}

\title[A 1/3-phenomenon for rhombus tilings]{A (conjectural) 1/3-phenomenon
for the number of rhombus
tilings of a hexagon which contain a fixed rhombus }
\author[C.~Krattenthaler]{\box\Adr}

\address{Institut f\"ur Mathematik der Universit\"at Wien,
Strudlhofgasse 4, A-1090 Wien, Austria.\newline
e-mail: KRATT@Ap.Univie.Ac.At\\
WWW: \tt http://www.mat.univie.ac.at/People/kratt}

\thanks {$^\dagger$ Research partially supported by the Austrian
Science Foundation FWF, grant P12094-MAT and P13190-MAT}

\begin{abstract}
We state, discuss, provide evidence for, and prove in special cases
the conjecture that the
probability that a random tiling by rhombi of a
hexagon with side lengths $2n+a,2n+b,2n+c,2n+a,2n+b,2n+c$ contains
the (horizontal) rhombus with coordinates $(2n+x,2n+y)$ is equal to $\frac {1}
{3}+g_{a,b,c,x,y}(n){\binom {2n}{n}}^3\Big/\binom {6n}{3n}$,
where $g_{a,b,c,x,y}(n)$ is a rational function in $n$. Several specific
instances of this ``$1/3$-phenomenon" are made explicit.
\end{abstract}

\maketitle

\begin{section}{Introduction and statement of the conjecture}
Let $a$, $b$ and $c$ be positive integers, and consider a
hexagon with side lengths $a,b,c,a,\break b,c$ whose angles are $120^\circ$
(see Figure~1.a for an example).
The subject of our interest is the enumeration of tilings
of this hexagon by rhombi (cf\@. Figure~1.b;
here, and in the sequel, by a
rhombus we always mean a rhombus with side lengths 1 and angles of
$60^\circ$ and $120^\circ$).

\begin{figure}[h]
\centertexdraw{
  \drawdim truecm  \linewd.02
  \rhombus \rhombus \rhombus \rhombus \ldreieck
  \move (-0.866025403784439 -.5)
  \rhombus \rhombus \rhombus \rhombus \rhombus \ldreieck
  \move (-1.732050807568877 -1)
  \rhombus \rhombus \rhombus \rhombus \rhombus \rhombus \ldreieck
  \move (-1.732050807568877 -1)
  \rdreieck
  \rhombus \rhombus \rhombus \rhombus \rhombus \rhombus \ldreieck
  \move (-1.732050807568877 -2)
  \rdreieck
  \rhombus \rhombus \rhombus \rhombus \rhombus \rhombus \ldreieck
  \move (-1.732050807568877 -3)
  \rdreieck
  \rhombus \rhombus \rhombus \rhombus \rhombus \rhombus
  \move (-1.732050807568877 -4)
  \rdreieck
  \rhombus \rhombus \rhombus \rhombus \rhombus
  \move (-1.732050807568877 -5)
  \rdreieck
  \rhombus \rhombus \rhombus \rhombus
\move(8 0)
\bsegment
  \drawdim truecm  \linewd.02
  \rhombus \rhombus \rhombus \rhombus \ldreieck
  \move (-0.866025403784439 -.5)
  \rhombus \rhombus \rhombus \rhombus \rhombus \ldreieck
  \move (-1.732050807568877 -1)
  \rhombus \rhombus \rhombus \rhombus \rhombus \rhombus \ldreieck
  \move (-1.732050807568877 -1)
  \rdreieck
  \rhombus \rhombus \rhombus \rhombus \rhombus \rhombus \ldreieck
  \move (-1.732050807568877 -2)
  \rdreieck
  \rhombus \rhombus \rhombus \rhombus \rhombus \rhombus \ldreieck
  \move (-1.732050807568877 -3)
  \rdreieck
  \rhombus \rhombus \rhombus \rhombus \rhombus \rhombus
  \move (-1.732050807568877 -4)
  \rdreieck
  \rhombus \rhombus \rhombus \rhombus \rhombus
  \move (-1.732050807568877 -5)
  \rdreieck
  \rhombus \rhombus \rhombus \rhombus
  \linewd.12
  \move(0 0)
  \RhombusA \RhombusB \RhombusB
  \RhombusA \RhombusA \RhombusB \RhombusA \RhombusB \RhombusB
  \move (-0.866025403784439 -.5)
  \RhombusA \RhombusB \RhombusB \RhombusB \RhombusB
  \RhombusA \RhombusA \RhombusB \RhombusA
  \move (-1.732050807568877 -1)
  \RhombusB \RhombusB \RhombusA \RhombusB \RhombusB \RhombusA
  \RhombusB \RhombusA \RhombusA
  \move (1.732050807568877 0)
  \RhombusC \RhombusC \RhombusC
  \move (1.732050807568877 -1)
  \RhombusC \RhombusC \RhombusC
  \move (3.464101615137755 -3)
  \RhombusC
  \move (-0.866025403784439 -.5)
  \RhombusC
  \move (-0.866025403784439 -1.5)
  \RhombusC
  \move (0.866025403784439 -2.5)
  \RhombusC \RhombusC
  \move (0.866025403784439 -3.5)
  \RhombusC \RhombusC \RhombusC
  \move (2.598076211353316 -5.5)
  \RhombusC
  \move (0.866025403784439 -5.5)
  \RhombusC
  \move (-1.732050807568877 -3)
  \RhombusC
  \move (-1.732050807568877 -4)
  \RhombusC
  \move (-1.732050807568877 -5)
  \RhombusC \RhombusC
\esegment
\htext (-1.5 -9){\small a. A hexagon with sides $a,b,c,a,b,c$,}
\htext (-1.5 -9.5){\small \hphantom{a. }where $a=3$, $b=4$, $c=5$}
\htext (6.8 -9){\small b. A rhombus tiling of a hexagon}
\htext (6.8 -9.5){\small \hphantom{b. }with sides $a,b,c,a,b,c$}
\rtext td:0 (4.3 -4){$\sideset {} c
    {\left.\vbox{\vskip2.2cm}\right\}}$}
\rtext td:60 (2.6 -.5){$\sideset {} {}
    {\left.\vbox{\vskip1.7cm}\right\}}$}
\rtext td:120 (-.44 -.25){$\sideset {}  {}
    {\left.\vbox{\vskip1.3cm}\right\}}$}
\rtext td:0 (-2.3 -3.6){$\sideset {c}  {}
    {\left\{\vbox{\vskip2.2cm}\right.}$}
\rtext td:240 (0 -7){$\sideset {}  {}
    {\left.\vbox{\vskip1.7cm}\right\}}$}
\rtext td:300 (3.03 -7.25){$\sideset {}  {}
    {\left.\vbox{\vskip1.4cm}\right\}}$}
\htext (-.9 0.1){$a$}
\htext (2.8 -.1){$b$}
\htext (3.2 -7.8){$a$}
\htext (-0.3 -7.65){$b$}
}
\caption{}
\end{figure}

As is well-known, the total number of rhombus tilings of a hexagon
with side lengths $a,b,c,a,b,c$ equals
$$\prod_{i=1}^a\prod_{j=1}^b\prod_{k=1}^c\frac{i+j+k-1}{i+j+k-2}.
$$
(This follows from MacMahon's enumeration
\cite[Sec.~429, $q\rightarrow 1$; proof in Sec.~494]{MacMAA}
of all plane partitions contained in an $a\times b\times c$ box,
as these are in bijection with rhombus tilings of a hexagon
with side lengths $a,b,c,a,b,c$, as explained e.g\@. in \cite{DaToAA}.)

The problem that we are going to address in this paper is the problem
of enumerating rhombus tilings of a hexagon which contain a given
fixed rhombus. Since the total number of rhombus tilings of a given
hexagon is known, thanks to MacMahon's formula, we may ask
equivalently the question of what the
probability is that a rhombus tiling of a hexagon that is chosen
uniformly at random (to be precise, it is the {\it tiling} which is chosen
at random, while the {\it hexagon} is given) contains a given fixed rhombus.
(For example, we may ask what the probability is that a randomly
chosen rhombus
tiling of the hexagon with side lengths $3,5,4,3,5,4$, shown in
Figure~4, contains the shaded rhombus. At this point the thick lines
are without relevance.)

\begin{figure}[h]
\centertexdraw{
  \drawdim truecm
  \move (4.33013 -2.5)
  \RhombusAsh
  \linewd.02
  \move (0 0)
  \rhombus \rhombus \rhombus \rhombus \rhombus \rhombus
  \rhombus \rhombus \rhombus \rhombus \rhombus
  \move (1.732050807568877 0)
  \rhombus \rhombus \rhombus \rhombus \rhombus
  \rhombus \rhombus \rhombus \rhombus
  \move (3.4641 0)
  \rhombus \rhombus \rhombus \rhombus
  \rhombus \rhombus \rhombus
  \move (5.19615 0)
  \rhombus \rhombus \rhombus
  \rhombus \rhombus
  \move (6.9282 0)
  \rhombus \rhombus
  \rhombus
  \move (8.66025 0)
  \rhombus
  \move (0 -1)
  \rhombus \rhombus \rhombus \rhombus \rhombus
  \rhombus \rhombus \rhombus \rhombus
  \move (0 -2)
  \rhombus \rhombus \rhombus \rhombus
  \rhombus \rhombus \rhombus
  \move (0 -3)
  \rhombus \rhombus \rhombus
  \rhombus \rhombus
  \move (0 -4)
  \rhombus \rhombus
  \rhombus
  \move (0 -5)
  \rhombus
  \move (0.866025403784439 .5)
  \StrichV
  \move (2.59808 .5)
  \StrichV
  \move (4.33013 .5)
  \StrichV
  \move (6.06218 .5)
  \StrichV
  \move (7.79423 .5)
  \StrichV
  \move (9.52628 .5)
  \StrichV
  \move (0.866025403784439 -5.5)
  \StrichA
  \move (2.59808 -5.5)
  \StrichA
  \move (4.33013 -5.5)
  \StrichA
  \move (6.06218 -5.5)
  \StrichA
  \move (7.79423 -5.5)
  \StrichA
  \move (9.52628 -5.5)
  \StrichA
  \move (0 0)
  \StrichX
  \move (0 -1)
  \StrichX
  \move (0 -2)
  \StrichX
  \move (0 -3)
  \StrichX
  \move (0 -4)
  \StrichX
  \move (0 -5)
  \StrichX
  \move (10.3923 0)
  \StrichY
  \move (10.3923 -1)
  \StrichY
  \move (10.3923 -2)
  \StrichY
  \move (10.3923 -3)
  \StrichY
  \move (10.3923 -4)
  \StrichY
  \move (10.3923 -5)
  \StrichY
}
\caption{}
\end{figure}

If this question is asked for an ``infinite" hexagon, i.e.,
if we imagine the 2-dimensional plane being covered by a triangular grid
(each triangle being an equilateral triangle; see Figure~2; at this
point shades in the figure should be ignored),
and ask the question of what the probability is that a particular
rhombus formed out of two adjacent triangles (for example the shaded
rhombus in Figure~2) is contained in a randomly chosen
rhombus tiling of the plane (that is compatible with the triangular
grid, of course), then there is a simple argument which shows that this
probability is $1/3$: Let us concentrate on one of the
two adjacent triangles out of which our fixed rhombus is formed. (In
Figure~3 we have enlarged the chosen rhombus. It is composed out of
the triangles labelled $0$ and $1$. We are going to concentrate on
the triangle labelled $0$.) This triangle is adjacent to exactly three
other triangles. (In Figure~3 these are the triangles labelled 1, 2
and 3.) In a rhombus tiling this triangle must be combined
with exactly one of these to form a rhombus in the tiling.
Hence, the probability that a {\it random} tiling will combine the
triangle with the particular one to obtain the fixed rhombus is $1/3$.

\begin{figure}[h]
\centertexdraw{
  \drawdim truecm \setunitscale2
  \linewd.02
  \RhombusAsh \ldreieck
  \move (0.866025403784439 .5)
  \ldreieck
  \move (0 0)
  \rhombus
\htext (.5 -.06){1}
\htext (1.15 -.06){0}
\htext (1.4 -.55){2}
\htext (1.4 .45){3}
}
\caption{}
\end{figure}

For a (finite) hexagon however, we must expect a very different
behaviour, resulting from the boundary of the hexagon.
The probability that a particular rhombus is contained in
a random tiling will heavily depend on where the rhombus is located
in the hexagon. (This is for example reflected in the asymptotic
result of Cohn, Larsen and Propp \cite[Theorem~1]{CoLPAA}.)
In particular, we must expect that the probability
will usually be different from $1/3$.

Rather surprisingly, John\footnote{In fact, in \cite{JohPAA} the
problem of finding the probability that, given a hexagonal graph, a
chosen fixed edge is contained in a randomly chosen perfect matching of the
graph is dealt with.
The motivation to consider this problem is that such hexagonal graphs
serve as models for benzenoid hydrocarbon molecules. The above probability is
called {\it Pauling's bond order}. It measures how stable a
carbon-carbon bond
(corresponding to the fixed edge) in a benzenoid hydrocarbon molecule
is.

It is well-known that this
problem is equivalent to our tiling problem. The link is a bijection
between rhombus tilings of a fixed subregion of the infinite triangular grid
(such as our hexagons) and perfect matchings of the hexagonal graph
which is, roughly speaking, the dual graph of the subregion (see e.g\@.
\cite{KupeAA}; ``roughly speaking" refers to the
little detail that the vertex corresponding to the outer face is
ignored in the dual graph construction).} \cite[bottom of p.~198]{JohPAA}
and Propp \cite[Problem~1]{PropAA,PropAH} made the empirical observation that in a hexagon
with side lengths $2n-1$, $2n-1$, $2n$, $2n-1$, $2n-1$, $2n$
the probability that the central rhombus is contained in a random
tiling is {\it exactly} $1/3$, the same being apparently true in a
hexagon with side lengths $2n$, $2n$, $2n-1$, $2n$, $2n$, $2n-1$.
These facts were proved by Ciucu and the
author \cite[Cor.~3]{CiKrAA} and, independently,
by Helfgott and Gessel \cite[Theorem~17]{HeGeAA}.
In fact, more generally, in both papers the probability that
in a hexagon with side lengths $N$, $N$, $M$, $N$, $N$, $M$,
$N\not\equiv M$ mod 2,
the central rhombus is contained in a random
tiling is expressed in terms of a single sum, from which the
$1/3$-result follows on simplification of the sum. These results were
generalized in two directions. On the one hand, Fulmek and the author
\cite{FuKrAC} found a single sum expression for this probability for
{\it any} rhombus on the
(horizontal) symmetry axis of the hexagon. On the other hand, Fischer
\cite{FiscAA} gave a single sum expression for the probability that
the central rhombus is contained in a random tiling of a hexagon with
{\it arbitrary} side lengths (i.e., with side lengths $a$, $b$, $c$,
$a$, $b$, $c$). Some further single sum expressions for probabilities
of ``near-central" rhombi to be contained in a random tiling of a
hexagon with sides $N$, $N$, $M$, $N$, $N$, $M$ have been derived in
\cite[Theorem~2]{FiscAA} and \cite{FuKrAD}.
Finally, in complete generality, Fischer
\cite[Lemma~2]{FiscAA} and Johansson \cite[(4.37)]{JohaAB} found {\it
triple sum} expressions for the
probability that a fixed (not necessarily central or near-central)
rhombus is contained in a random
tiling of a hexagon with side lengths $a$, $b$, $c$, $a$, $b$, $c$.
(These two triple sum expressions are completely different from each
other.)

The purpose of this paper is to report a curious manifestation of the
fact that ``in the limit" the probability that a particular rhombus
is contained in a random tiling is $1/3$. Roughly speaking, it seems
that the probability equals
$$1/3\text { plus a ``nice" expression.}$$
To make
this precise, we need to introduce a convention of how to describe
the position of a rhombus in
a given hexagon. First of all, without loss of generality, we may
restrict our considerations to the case where the fixed rhombus is a
horizontal rhombus (by which we mean a rhombus such as the shaded
ones in Figures~2--4), which we shall do
for the rest of the paper. (The other two
types of rhombi are then covered via a rotation by $120^\circ$, respectively
by $240^\circ$.) In order to describe the position of a rhombus in
the hexagon, we introduce, following \cite{FiscAA},
the following oblique angled coordinate
system: Its origin is located in one of the two vertices where the
sides of lengths $b$ and $c$ meet, and the axes are induced by those
two sides (see Figure~4). The units are chosen such that the
grid points of the triangular grid are exactly the integer points
in this coordinate system. (That is to say, the
two triangles in Figure~4 with vertices in the origin form the unit
`square.') Thus, in this coordinate system, the bottom-most point of
the shaded hexagon in Figure~4
has coordinates $(5,4)$.
\newbox\obliquebox
\setbox\obliquebox\hbox{\small
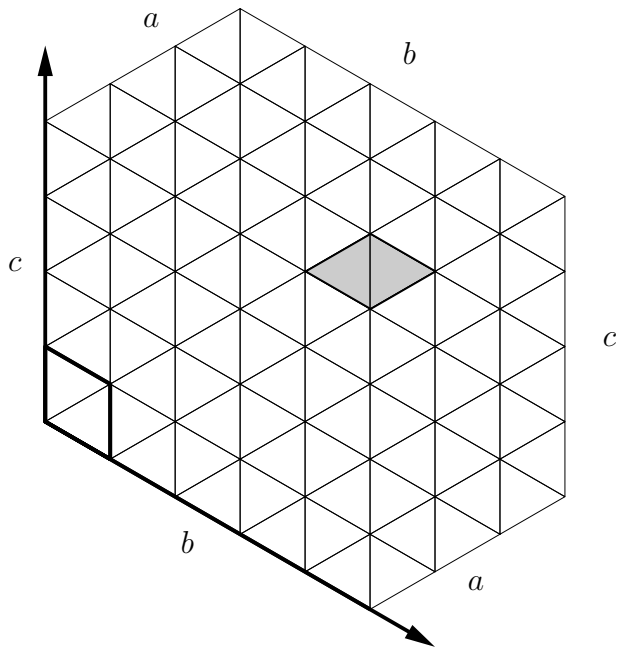
The oblique angled coordinate system}
\begin{figure}
\centertexdraw{
\drawdim truecm
\arrowheadtype t:F
\move (3.4641 -1)
\RhombusAsh
\move (0 0)
\linewd 0.05
\move(0 -3)
\ravec(0 5)
\move(0 -3)
\rlvec(0.8660254037844 -0.5)
\rlvec(0 1)
\rlvec(-0.866025403784 0.5)
\rlvec(0 -1)
\rlvec(0.8660254037844 -0.5)
\rlvec(0.8660254037844 -0.5)
\rlvec(0.8660254037844 -0.5)
\rlvec(0.8660254037844 -0.5)
\rlvec(0.8660254037844 -0.5)
\ravec(0.8660254037844 -0.5)
\linewd 0.01
\move(0 0)
\DreieckBreit
\rmove(0 1)
\DreieckSpitz
\rmove(0.8660254037844 -0.5)
\DreieckBreit
\rmove(0 1)
\DreieckSpitz
\rmove(0.8660254037844 -0.5)
\DreieckBreit
\rmove(0 1)
\DreieckSpitz
\rmove(0.8660254037844 -0.5)
\DreieckBreit
\rmove(-0.8660254037844 -0.5)
\rmove(-0.8660254037844 -0.5)
\rmove(-0.8660254037844 -0.5)
\rmove(0 -1)
\DreieckBreit
\rmove(0 1)
\DreieckSpitz
\rmove(0.8660254037844 -0.5)
\DreieckBreit
\rmove(0 1)
\DreieckSpitz
\rmove(0.8660254037844 -0.5)
\DreieckBreit
\rmove(0 1)
\DreieckSpitz
\rmove(0.8660254037844 -0.5)
\DreieckBreit
\rmove(0 1)
\DreieckSpitz
\rmove(0.8660254037844 -0.5)
\DreieckBreit
\rmove(-0.8660254037844 -0.5)
\rmove(-0.8660254037844 -0.5)
\rmove(-0.8660254037844 -0.5)
\rmove(-0.8660254037844 -0.5)
\rmove(0 -1)
\DreieckBreit
\rmove(0 1)
\DreieckSpitz
\rmove(0.8660254037844 -0.5)
\DreieckBreit
\rmove(0 1)
\DreieckSpitz
\rmove(0.8660254037844 -0.5)
\DreieckBreit
\rmove(0 1)
\DreieckSpitz
\rmove(0.8660254037844 -0.5)
\DreieckBreit
\rmove(0 1)
\DreieckSpitz
\rmove(0.8660254037844 -0.5)
\DreieckBreit
\rmove(0 1)
\DreieckSpitz
\rmove(0.8660254037844 -0.5)
\DreieckBreit
\rmove(-0.8660254037844 -0.5)
\rmove(-0.8660254037844 -0.5)
\rmove(-0.8660254037844 -0.5)
\rmove(-0.8660254037844 -0.5)
\rmove(-0.8660254037844 -0.5)
\rmove(0 -1)
\DreieckBreit
\rmove(0 1)
\DreieckSpitz
\rmove(0.8660254037844 -0.5)
\DreieckBreit
\rmove(0 1)
\DreieckSpitz
\rmove(0.8660254037844 -0.5)
\DreieckBreit
\rmove(0 1)
\DreieckSpitz
\rmove(0.8660254037844 -0.5)
\DreieckBreit
\rmove(0 1)
\DreieckSpitz
\rmove(0.8660254037844 -0.5)
\DreieckBreit
\rmove(0 1)
\DreieckSpitz
\rmove(0.8660254037844 -0.5)
\DreieckBreit
\rmove(0 1)
\DreieckSpitz
\rmove(0.8660254037844 -0.5)
\DreieckBreit
\rmove(-0.8660254037844 -0.5)
\rmove(-0.8660254037844 -0.5)
\rmove(-0.8660254037844 -0.5)
\rmove(-0.8660254037844 -0.5)
\rmove(-0.8660254037844 -0.5)
\rmove(-0.8660254037844 -0.5)
\DreieckSpitz
\rmove(0.8660254037844 -0.5)
\DreieckBreit
\rmove(0 1)
\DreieckSpitz
\rmove(0.8660254037844 -0.5)
\DreieckBreit
\rmove(0 1)
\DreieckSpitz
\rmove(0.8660254037844 -0.5)
\DreieckBreit
\rmove(0 1)
\DreieckSpitz
\rmove(0.8660254037844 -0.5)
\DreieckBreit
\rmove(0 1)
\DreieckSpitz
\rmove(0.8660254037844 -0.5)
\DreieckBreit
\rmove(0 1)
\DreieckSpitz
\rmove(0.8660254037844 -0.5)
\DreieckBreit
\rmove(0 1)
\DreieckSpitz
\rmove(0.8660254037844 -0.5)
\DreieckBreit
\rmove(-0.8660254037844 -0.5)
\rmove(-0.8660254037844 -0.5)
\rmove(-0.8660254037844 -0.5)
\rmove(-0.8660254037844 -0.5)
\rmove(-0.8660254037844 -0.5)
\rmove(-0.8660254037844 -0.5)
\DreieckSpitz
\rmove(0.8660254037844 -0.5)
\DreieckBreit
\rmove(0 1)
\DreieckSpitz
\rmove(0.8660254037844 -0.5)
\DreieckBreit
\rmove(0 1)
\DreieckSpitz
\rmove(0.8660254037844 -0.5)
\DreieckBreit
\rmove(0 1)
\DreieckSpitz
\rmove(0.8660254037844 -0.5)
\DreieckBreit
\rmove(0 1)
\DreieckSpitz
\rmove(0.8660254037844 -0.5)
\DreieckBreit
\rmove(0 1)
\DreieckSpitz
\rmove(0.8660254037844 -0.5)
\DreieckBreit
\rmove(0 1)
\DreieckSpitz
\rmove(0.8660254037844 -0.5)
\rmove(-0.8660254037844 -0.5)
\rmove(-0.8660254037844 -0.5)
\rmove(-0.8660254037844 -0.5)
\rmove(-0.8660254037844 -0.5)
\rmove(-0.8660254037844 -0.5)
\rmove(-0.8660254037844 -0.5)
\DreieckSpitz
\rmove(0.8660254037844 -0.5)
\DreieckBreit
\rmove(0 1)
\DreieckSpitz
\rmove(0.8660254037844 -0.5)
\DreieckBreit
\rmove(0 1)
\DreieckSpitz
\rmove(0.8660254037844 -0.5)
\DreieckBreit
\rmove(0 1)
\DreieckSpitz
\rmove(0.8660254037844 -0.5)
\DreieckBreit
\rmove(0 1)
\DreieckSpitz
\rmove(0.8660254037844 -0.5)
\DreieckBreit
\rmove(0 1)
\DreieckSpitz
\rmove(0.8660254037844 -0.5)
\rmove(-0.8660254037844 -0.5)
\rmove(-0.8660254037844 -0.5)
\rmove(-0.8660254037844 -0.5)
\rmove(-0.8660254037844 -0.5)
\rmove(-0.8660254037844 -0.5)
\DreieckSpitz
\rmove(0.8660254037844 -0.5)
\DreieckBreit
\rmove(0 1)
\DreieckSpitz
\rmove(0.8660254037844 -0.5)
\DreieckBreit
\rmove(0 1)
\DreieckSpitz
\rmove(0.8660254037844 -0.5)
\DreieckBreit
\rmove(0 1)
\DreieckSpitz
\rmove(0.8660254037844 -0.5)
\DreieckBreit
\rmove(0 1)
\DreieckSpitz
\rmove(0.8660254037844 -0.5)
\rmove(-0.8660254037844 -0.5)
\rmove(-0.8660254037844 -0.5)
\rmove(-0.8660254037844 -0.5)
\rmove(-0.8660254037844 -0.5)
\DreieckSpitz
\rmove(0.8660254037844 -0.5)
\DreieckBreit
\rmove(0 1)
\DreieckSpitz
\rmove(0.8660254037844 -0.5)
\DreieckBreit
\rmove(0 1)
\DreieckSpitz
\rmove(0.8660254037844 -0.5)
\DreieckBreit
\rmove(0 1)
\DreieckSpitz
\rmove(0.8660254037844 -0.5)
\move(-0.5 -1)
\htext{$c$}
\move(1.299038105677 1.75)
\rmove(0  0.5)
\htext{$a$}
\rmove(0 -0.5)
\rmove(1.299038105677 0.75)
\rmove(1.299038105677 -0.75)
\rmove(0.8660254037844 -0.5)
\rmove(0 0.5)
\htext{$b$}
\rmove(0 -0.5)
\rmove(1.299038105677 -0.75)
\rmove(0.8660254037844 -0.5)
\rmove(0 -2)
\rmove(0.5 0)
\htext{$c$}
\rmove(-0.5 0)
\rmove(0 -2)
\rmove(-1.299038105677 -0.75)
\rmove(0 -0.5)
\htext{$a$}
\rmove(0.5 0)
\rmove(-1.299038105677 -0.75)
\rmove(-1.299038105677 0.75)
\rmove(-0.8660254037844 0.5)
\rmove(-0.8660254037844 0.5)
\rmove(0 -0.5)
\htext{$b$}
}
\caption{\unhbox\obliquebox}
\end{figure}

With this convention, we have the following conjecture. It extends an
(ex)conj\-ect\-ure by Propp \cite[Problem~4]{PropAA,PropAH}.
\begin{Conjecture} Let $a$, $b$, $c$, $x$ and $y$ be arbitrary integers.
Then the probability that a randomly chosen rhombus tiling of a
hexagon with side lengths $2n+a$,
$2n+b$, $2n+c$, $2n+a$, $2n+b$, $2n+c$ contains
the (horizontal) rhombus with bottom-most vertex $(2n+x,2n+y)$ (in
the oblique angled coordinate system)
is equal to
\begin{equation} \label{e1.1}
\frac {1}
{3}+f_{a,b,c,x,y}(n){\binom {2n}{n}}^3\bigg/\binom {6n+2}{3n+1}
\quad \quad \text {for $n>n_0$,}
\end{equation}
for a suitable $n_0$ which depends on $a$, $b$, $c$, $x$ and $y$,
where $f_{a,b,c,x,y}(n)$ is a rational function in
$n$.\footnote{\rm This statement is clearly equivalent to the statement
in the
abstract. The form \eqref{e1.1} of the expression is more convenient in the
subsequent listing of special cases.}
\end{Conjecture}

As is shown in Section~2, for any {\it specific} $a$, $b$, $c$, $x$, $y$
the corresponding formula for $f_{a,b,c,x,y}(n)$ can be worked out
completely automatically by the use of a computer
(given that the Conjecture is true, of course).
We have in fact produced a huge
list of such formulas, of which we list a few selected instances below.
As we explain in Section~3, any of these is (at least) a
``near-theorem," in the sense that it could be proved automatically
by the available multisum algorithms, provided
there is enough computer memory available (and, thus, will at least
be a theorem in the near future). Also in Section~3, we elaborate more
precisely on which of these are just conjectural, and which
of them are already theorems\footnote{For the convenience of the
reader, we have marked conjectures by an asterisk in the equation
number.}. However, we do not know how to prove
the Conjecture {\it in general\/}, that is, for {\it generic} values
of $a$, $b$, $c$, $x$, and $y$ (cf\@. Section~3 for a possible
approach).

Here is the announced excerpt from our list of special instances of
the Conjecture:

\begin{equation} \label{L1.1}
f_{-1, -1, 0, -1, -1}(n) = f_{2, 2, 1, 2, 1}(n) = 0\quad \quad
\text{for $n\ge 1$},
\end{equation}
\begin{multline} \label{L1.2}
f_{2, 1, 1, 2, 1}(n) = f_{2, 1, 1, 1, 1}(n) =
  f_{1, 2, 1, 2, 1}(n) = f_{1, 2, 1, 1, 0}(n) \\
=
  f_{-1, 0, 0, 0, -1}(n) = f_{-1, 0, 0, -1, -1}(n) =
  f_{0, -1, 0, 0, 0}(n) = f_{0, -1, 0, -1, -1}(n) = 0\quad \quad
\text{for $n\ge 1$},
\end{multline}
\begin{multline} \label{L1.3}
f_{1, 1, 1, 1, 1}(n) =
f_{1, 1, 0, 1, -1}(n) =
f_{1, 1, 0, 1, 1}(n)\\
 =
f_{1, 1, 1, 0, -1}(n) =
f_{1, 1, 1, 0, 1}(n) =
f_{1, 1, 1, 1, 0}(n) =
  f_{0, 2, 0, 1, 0}(n) \kern.8cm\\
=
f_{2, 0, 0, 0, 1}(n) =
f_{2, 0, 0, 1, 0}(n) =
f_{2, 0, 0, 1, 1}(n) =
f_{2, 0, 0, 2, 0}(n) =
 \frac 1 3\quad \quad \text{for $n\ge 1$},
\end{multline}
\def\theequation{\mbox{\thesection.\arabic{equation}*}}
\begin{equation} \label{L1.3a}
f_{1, 1, 0, 0, -1}(n) =
f_{1, 1, 0, 0, 0}(n) =
  f_{2, 0, 1, 1, 0}(n) =
 \frac 1 3\quad \quad \text{for $n\ge 1$},
\end{equation}
\def\theequation{\mbox{\thesection.\arabic{equation}}}
\begin{equation} \label{L1.4}
f_{0, 2, 1, 1, 0}(n) =
-  \frac {2} 3\quad \quad \text{for $n\ge 1$},
\end{equation}
\def\theequation{\mbox{\thesection.\arabic{equation}*}}
\begin{equation} \label{L1.4a}
f_{1, 1, 1, 0, 0}(n) =
-  \frac {2} 3\quad \quad \text{for $n\ge 1$},
\end{equation}
\def\theequation{\mbox{\thesection.\arabic{equation}}}
\begin{equation} \label{L1.5}
f_{1, 1, 0, 1, 0}(n) =
  \frac 4 3\quad \quad \text{for $n\ge 1$},
\end{equation}
\def\theequation{\mbox{\thesection.\arabic{equation}*}}
\begin{equation} \label{L1.5a}
f_{4, 3, 1, 3, 2}(n) =
f_{4, 3, 1, 4, 2}(n) =
  \frac 4 3\quad \quad \text{for $n\ge 1$},
\end{equation}
\def\theequation{\mbox{\thesection.\arabic{equation}}}
\begin{equation} \label{L1.6a}
f_{0, 0, 0, 0, 0}(n) =
-  \frac { (6 n + 1)} {
   6 (3 n + 1)}\quad \quad \text{for $n\ge 1$},
\end{equation}
\begin{equation} \label{L1.6}
f_{0, 0, 1, 0, 0}(n) =
-  \frac {2 (6 n + 1)} {
   3 (3 n + 1)}\quad \quad \text{for $n\ge 1$},
\end{equation}
\begin{equation} \label{L1.7}
f_{3, 3, 0, 3, 1}(n) =
  \frac {2 (2 n + 1) (3 n + 2) (4 n + 5)} {
   3 (n + 1)^2 (6 n + 5)}\quad \quad \text{for $n\ge 1$},
\end{equation}
\def\theequation{\mbox{\thesection.\arabic{equation}*}}
\begin{equation} \label{L1.8}
f_{2, 1, 0, 3, -1}(n) =
  \frac {4 n^3+ 18 n^2 + 12 n +1} {6 (n+1) ^2 (2 n-1)}\quad \quad
\text{for $n\ge 2$},
\end{equation}
\begin{equation} \label{L1.9}
f_{5, 1, 0, 3, 2}(n) =
  \frac {(3 n + 2) (16 n^3 + 54 n^2 + 57 n + 20)} {
   3 (n + 1)^2 (n + 2) (6 n + 5)}\quad \quad \text{for $n\ge 1$},
\end{equation}
\begin{equation} \label{L1.10}
f_{-1, 5, 0, 2, -1}(n) =
  \frac {(3 n + 2) (2 n^2 + 4 n + 1)} {
   3 (n + 1)^2 (n + 2)}\quad \quad \text{for $n\ge 1$},
\end{equation}
\begin{multline} \label{L1.11}
f_{10,3,0,1,4}(n) =
\tfrac {(2 n + 1) (2 n + 3) (3 n + 2) (3 n + 4) (3 n + 5)
  } {
  6 (n + 1)^2 (n + 2)^2 (n + 3)^2 (n + 4)^2 (n + 5) (2 n -3)
    (2 n -1) (6 n + 5) (6 n + 7) (6 n + 11)}\\
\scriptstyle\times
(176 n^9 + 3080 n^8 + 21692 n^7 + 74546 n^6 +
      102578 n^5  - 73279 n^4  - 362598 n^3 - 283977 n^2 + 24762 n + 55440
)\\\text{for $n\ge 2$}.
\end{multline}
\def\theequation{\mbox{\thesection.\arabic{equation}}}

\end{section}

\begin{section}{How are these conjectures and results discovered?}
Point of departure for all these discoveries is an observation by
Propp \cite[Problem~4]{PropAA,PropAH}: He conjectured that the
probability that a randomly chosen rhombus tiling of a hexagon with {\it
all\/} side lengths equal to $N$
contains the ``near-central" rhombus (this is the rhombus with
bottom-most vertex $(N,N)$ in the oblique angled coordinate system)
is equal to $1/3$ plus a ``nice"
formula in $N$.\footnote{Commonly, by a ``nice" formula
one means an expression which is
built by forming products and quotients of factorials. A strong
indication that one encounters a sequence $(a_N)_{N\ge0}$
for which a ``nice" formula
exists is that the prime factors in the prime factorization of $a_N$ do not
grow rapidly as $N$ becomes larger. (In fact, they should grow
linearly.)} Should this observation be true, then the {\sl
Mathematica} program {\tt Rate}\footnote{written by the author;
available from {\tt
http://radon.mat.univie.ac.at/People/kratt}; see \cite[Appendix~A]{KratBN}
for an explanation of how the program works.} (``Rate!" is
German for ``Guess!"), respectively its {\sl Maple} equivalent {\tt
GUESS}\footnote{written by
Fran\c cois B\'eraud and Bruno Gauthier; available from
{\tt http://www-igm.univ-mlv.fr/\~{}gauthier}.},
will find the formula, given enough initial terms of the sequence.

Let us see how this
works in the case that $N$ is odd.
For the generation of the probabilities, Propp used the
programs {\tt vaxmaple}\footnote{written by Greg Kuperberg, Jim Propp
and David Wilson; available at {\tt
http://math.wisc.edu/\~{}propp/software.html}.} and
{\tt vaxmacs}\footnote{written by
David Wilson; also available at {\tt
http://math.wisc.edu/\~{}propp/software.html}.},
which are based on the evaluation
of determinants of large (if though sparse) matrices.
However, since then triple sum formulas have been found by
Fischer \cite[Lemma~2]{FiscAA} and Johansson \cite[(4.37)]{JohaAB},
which allow to generate these probabilities much more efficiently.
We choose to use Fischer's formula. We state it below.

\begin{Theorem}
Let $a$, $b$ and $c$ be positive integers, and
let $(x,y)$ be an integer point such that $0\le x\le b+a-1$ and $1\le
y\le c+a-1$. Then the probability that a randomly chosen rhombus tiling of a
hexagon with side lengths $a$, $b$, $c$, $a$, $b$, $c$ contains
the (horizontal) rhombus with bottom-most vertex $(x,y)$ (in
the oblique angled coordinate system)
is equal to
\begin{multline} \label{e2.1}
{\frac {c!} {({ \textstyle b+1}) _{c} }}
\sum _{i=1} ^{a} \sum _{j=1} ^{a} \sum _{s=1} ^{j}
{{\left( -1 \right)      }^{i + s}}
{\binom {j-1} {s-1}}\,{\binom {c + i + x - y-2} {x-1
      }}\,{\binom {b + s - x + y-1} {b + s - x-1}}\\
\cdot
{\frac {  ({ \textstyle b +
      1}) _{s-1} \,({ \textstyle c+1}) _{ i-1} \,({ \textstyle b + c +
      i}) _{j - i} } {\left( j - i \right) !\,\left( i-1 \right) !\,({
      \textstyle b + c+1}) _{s-1} }}.
\end{multline}
\end{Theorem}

We now program this formula in {\sl Mathematica}.
\MATH
\goodbreakpoint%
Mathematica 2.2 for DOS 387
Copyright 1988-93 Wolfram Research, Inc.
\goodbreakpoint%
In[1]:= F[a\MATHtief ,b\MATHtief ,c\MATHtief ,x\MATHtief ,y\MATHtief ]:=c!/Pochhammer[b+1,c]*
\leavevmode   Sum[Sum[Sum[(-1)\MATHhoch (i+s)*Binomial[j-1,s-1]*
\leavevmode     Binomial[c+i+x-y-2,x-1]*Binomial[b+s-x+y-1,b+s-x-1]*
\leavevmode     Pochhammer[b+1,s-1]*Pochhammer[c+1,i-1]*
\leavevmode     Pochhammer[b+c+i,j-i]/(j-i)!/(i-1)!/Pochhammer[b+c+1,s-1],
\leavevmode   %
\MATHlbrace s,1,j%
\MATHrbrace ],%
\MATHlbrace j,1,a%
\MATHrbrace ],%
\MATHlbrace i,1,a%
\MATHrbrace ]
\goodbreakpoint%
\endMATH
\noindent
Now we generate the first eleven values of these probabilities for
$N=2n+1$ and subtract $1/3$ from them.
\MATH
\goodbreakpoint%
In[2]:= Table[F[2n+1,2n+1,2n+1,2n+1,2n+1]-1/3,%
\MATHlbrace n,1,11%
\MATHrbrace ]
\goodbreakpoint%
          4    3    2000    245    296352    142296    43188288
Out[2]= %
\MATHlbrace ---, ---, ------, -----, --------, --------, ----------,
\leavevmode         105  143  138567  22287  33393355  19126225  6743906935

\leavevmode      759169125     15365378600    55469016746     805693639296
\MATHgroesser     ------------, -------------, --------------, ---------------%
\MATHrbrace
\leavevmode     135054066707  3067656658059  12280863528759  195909013434965
\goodbreakpoint%
\endMATH
\noindent
Next we load {\tt Rate}, and apply {\tt Rate}'s function {\tt Ratekurz}
to the sequence of numbers.
\MATH
\goodbreakpoint%
In[3]:= \MATHkleiner \MATHkleiner rate.m
\goodbreakpoint%
In[4]:= Apply[Ratekurz,\%2]
\goodbreakpoint%
\leavevmode                             2
\leavevmode                   (1 + 2 i1)  (2 + 3 i1) (4 + 3 i1)
\leavevmode         4 Product[---------------------------------,\,%
\MATHlbrace i1, 1, -1 + i0%
\MATHrbrace ]
\leavevmode                            2
\leavevmode                    (1 + i1)  (5 + 6 i1) (7 + 6 i1)
Out[4]= %
\MATHlbrace -------------------------------------------------------------%
\MATHrbrace
\leavevmode                                      105
\goodbreakpoint%
\endMATH
\noindent
The program outputs a formula which generates the terms of the
sequence that was given as an input. The formula is written
as a function in {\tt i0}, i.e., we
must replace {\tt i0} by $n$. In more compact terms, the formula can
be rewritten as
\begin{equation} \label{e2.2}
\frac {1} {3}{\binom {2n}{n}}^3\bigg/\binom {6n+2}{3n+1}.
\end{equation}
(It should be observed that this expression is exactly
the one which features in \eqref{e1.1}.)

At this point, this formula is of course just a conjecture. It has
however been proved in \cite[Corollary~7, (1.9)]{FuKrAD}.

Being adventurous, one tries the same thing for other choices of the
parameters $a$, $b$, $c$, $x$ and $y$.
Very quickly one discovers, that a similar phenomenon
seems to occur for any choice $2n+a$, $2n+b$, $2n+c$ for the side
lengths and $(2n+x,2n+y)$ for the coordinates of the bottom-most
point of the fixed rhombus, where $a$, $b$, $c$, $x$, and $y$ are fixed
integers. Although the (conjectural)
expressions that one finds need not be ``nice" anymore in the
strict sense above, it is at worst polynomial factors in $n$ that
appear in addition. Moreover, one also realizes soon that division of
such an expression by the expression in \eqref{e2.2} apparently always
results in a rational function in $n$, i.e., the Conjecture in
Section~1 is discovered.

Let us see just one such example. We choose a hexagon with side
lengths $2n+2$, $2n+1$, $2n$, $2n+2$, $2n+1$, $2n$, and
$(2n+3,2n-1)$ for the coordinates of the bottom-most point of the
fixed rhombus. Then we obtain the following numbers for
$n=1,2,\dots,15$. The reader should note that we immediately divide
the expression \eqref{e2.2}.
\MATH
\goodbreakpoint%
In[5]:= Table[(F2[2n+2,2n+1,2n,2n+3,2n-1]-1/3)/(Binomial[2n,n]\MATHhoch 3/
\leavevmode        Binomial[6n+2,3n+1]),%
\MATHlbrace n,1,15%
\MATHrbrace ]
\goodbreakpoint%
\leavevmode           35   43  307  593   337  1585  2339  1099  4483   5921
Out[5]= %
\MATHlbrace -(--), --, ---, ----, ---, ----, ----, ----, -----, -----,
\leavevmode           12   54  480  1050  648  3234  4992  2430  10200  13794

\leavevmode     2545  9649   11987  4891   17731
\MATHgroesser     ----, -----, -----, -----, -----%
\MATHrbrace
\leavevmode     6048  23322  29400  12150  44544
\goodbreakpoint%
\endMATH
\noindent
By having a brief glance at this sequence, it seems that the first
term is ``alien," so let us better drop it.
\MATH
\goodbreakpoint%
In[6]:= Drop[\%,1]
\goodbreakpoint%
\leavevmode         43  307  593   337  1585  2339  1099  4483   5921
Out[6]= %
\MATHlbrace --, ---, ----, ---, ----, ----, ----, -----, -----,
\leavevmode         54  480  1050  648  3234  4992  2430  10200  13794

\leavevmode     2545  9649   11987  4891   17731
\MATHgroesser     ----, -----, -----, -----, -----%
\MATHrbrace
\leavevmode     6048  23322  29400  12150  44544
\goodbreakpoint%
\endMATH
\noindent
By the discussion above, this should be a sequence which is given by
a rational function in $n$. Therefore is suffices to apply {\tt Rate}'s {\tt
Rateint} (which does just rational interpolation, in contrast to {\tt
Ratekurz}, which tries several other things, and which is therefore slower).
\MATH
\goodbreakpoint%
In[7]:= Apply[Rateint,\%]
\goodbreakpoint%
\leavevmode                           2       3
\leavevmode         35 + 60 i0 + 30 i0  + 4 i0
Out[7]= %
\MATHlbrace ---------------------------%
\MATHrbrace
\leavevmode                     2
\leavevmode           6 (2 + i0)  (1 + 2 i0)
\goodbreakpoint%
\endMATH
Again, the program outputs the formula as a function in {\tt i0}.
Since initially we dropped the first term of the sequence, we must
now replace {\tt i0} by $n-1$.
\MATH
\goodbreakpoint%
In[8]:= Factor[\%/.i0-\MATHgroesser n-1]
\goodbreakpoint%
\leavevmode                        2      3
\leavevmode         1 + 12 n + 18 n  + 4 n
Out[8]= %
\MATHlbrace -----------------------%
\MATHrbrace
\leavevmode                   2
\leavevmode          6 (1 + n)  (-1 + 2 n)
\goodbreakpoint%
\endMATH
Hence, if the Conjecture in Section~1 is true, $f_{2,1,0,3,-1}(n)$ must
be the expression given in the output {\tt Out[8]}. (Thus, we have
discovered Eq.~(1.13).) Again, at this
point, this is just a conjecture.

The Equations
\eqref{L1.1}--(1.16) in Section~1 are all found in the same way.
\end{section}

\begin{section}{Discussion: How to prove the conjecture?}
A possible approach to prove the Conjecture in Section~1 is to start
with the expression \eqref{e2.1} (or with the alternative expression
\cite[(4.37)]{JohaAB}), replace $a$ by $2n+a$,
$b$ by $2n+b$, $c$ by $2n+c$, $x$ by $2n+x$, $y$ by $2n+y$,
and by some manipulation (for example, by applying
hypergeometric transformation and summation formulas) convert it into
the form \eqref{e2.2}. Everybody who has some experience with manipulating
binomial/hypergeometric sums will immediately realize that this is a
formidable task. In particular, it seems a bit mysterious how one
should be able to isolate ``$1/3$" from the ``rest." In any case, I do
not know how to prove the Conjecture in this manner, nor in any other
way.

On the other hand, as we explained in Section~2, for any {\it
specific}
values of $a$, $b$, $c$, $x$, and $y$, it is routine to find a
conjectural expression for the rational function $f_{a,b,c,x,y}(n)$
(given that the Conjecture is true). In turn, once such an expression
is available, it can (at least in principle) be verified completely
automatically. For, what one has to prove is the equality of the
expressions \eqref{e2.1}, with the above replacements, and \eqref{e1.1}, where
$f_{a,b,c,x,y}(n)$ is the explicit rational function found by the
computer. That is to say, one has to prove that a certain triple sum
equals a closed form expression. Clearly, this can be done (again,
at least in principle) by the available multisum
algorithms\footnote{The first (theoretical) algorithm for proving
multisum identities automatically was given by Wilf and Zeilberger
\cite{WiZeAC}. A considerable enhancement and speedup was
accomplished by
Wegschaider \cite{WegsAA}, who combined the ideas of Wilf and
Zeilberger with ideas of Verbaeten \cite{VerbAA}. Wegschaider's {\sl
Mathematica} implementation is available from {\tt
http://www.risc.uni-linz.ac.at/research/combinat/risc/software}.},
by using the algorithm to find a recurrence in $n$ for the expression
\eqref{e2.1}, and subsequently checking that the expression \eqref{e1.1}, with the
computer guess for $f_{a,b,c,x,y}(n)$, satisfies the same recurrence.
Unfortunately, in any case that I tried, the computer ran out of
memory.

However, as we already mentioned in the Introduction, in some cases
formulas in form of single hypergeometric sums are available. If one
is in such a case then one would proceed as in the above paragraph,
but one would replace the multisum algorithm by Zeilberger's
algorithm\footnote{A {\it Maple}
implementation written by Doron
Zeilberger is available from {\tt
http://www.math.temple.edu/\~{}zeilberg}; a {\sl Mathematica}
implementation
written by Markus Schorn and Peter Paule is available from {\tt
http://www.risc.uni-linz.ac.at/research/combinat/risc/software}.}
(see \cite{PaScAA,PeWZAA,ZeilAM,ZeilAV}). The advantage is
that, in contrast to the multisum algorithm, Zeilberger's algorithm is
very efficient. At any rate, in any case that I looked at in
connection with our problem, the Zeilberger algorithm was
successful. That is to say, if I am allowed to somewhat overstate it,
{\it whenever one is in a
case where a single sum formula is available, one has a theorem}
(i.e., Zeilberger's algorithm {\it will\/} prove
that the empirical found rational function $f_{a,b,c,x,y}(n)$
does indeed satisfy the Conjecture for {\it all\/} values of $n$).

For the sake of completeness, we list the vectors $(a,b,c,x,y)$
for which single sums are available for $f_{a,b,c,x,y}(n)$. Clearly,
it suffices to restrict $c$ to $0$ and $1$. (All other values can be
attained by shifts of $n$.)
\begin{enumerate}
\item[(A)] By \cite[Theorem~1, (1.2)]{FiscAA}:
$(2a'+1,2b'+1,0,a'+b'+1,a')$, for integers $a'$ and
$b'$.
\item[(B)] By \cite[Theorem~1, (1.3)]{FiscAA}:
$(2a',2b',1,a'+b',a')$, for integers $a'$ and
$b'$.
\item[(C)] By \cite[Theorem~2, (1.4)]{FiscAA}:
$(2a'+1,2b'+1,1,a'+b'+1,a'+1)$ and
$(2a'+1,2b'+1,1,a'+b'+1,a')$, for integers $a'$ and
$b'$.
\item[(D)] By \cite[Theorem~2, (1.5)]{FiscAA}:
$(2a',2b',0,a'+b',a')$ and
$(2a',2b',0,a'+b',a'-1)$, for integers $a'$ and
$b'$.
\item[(E)] By \cite[Theorem~1]{FuKrAC}:
$(a',a',0,2x'+1,x')$, for integers $a'$ and
$x'$.
\item[(F)] By \cite[Theorem~2]{FuKrAC}:
$(a',a',1,2x',x'-1)$, for integers $a'$ and
$x'$.
\item[(G)] By \cite[Theorem~3]{FuKrAD}:
$(2a',2a',1,2a',a'+1)$ and
$(2a',2a',1,2a',a'-1)$, for an integer $a'$.
\item[(H)] By \cite[Theorem~4]{FuKrAD}:
$(2a'+1,2a'+1,0,2a'+1,a'+1)$ and
$(2a'+1,2a'+1,0,2a'+1,a'-1)$, for an integer $a'$.
\item[(I)] By \cite[Theorem~5]{FuKrAD}:
$(2a',2a',0,2a',a'+1)$ and
$(2a',2a',0,2a',a'-2)$, for an integer $a'$.
\item[(J)] By \cite[Theorem~6]{FuKrAD}:
$(2a'+1,2a'+1,1,2a'+1,a'+2)$ and
$(2a'+1,2a'+1,1,2a'+1,a'-1)$, for an integer $a'$.
\end{enumerate}

Thus, choosing $a'=b'=0$ in (A), we see for example that the
expression for $f_{0,0,1,0,0}(n)$ given in \eqref{L1.6} is in fact a
theorem. For, by Theorem~1, (1.3) in \cite{FiscAA} with $a=b=2n$,
$c=2n+1$, the
probability that a randomly chosen rhombus tiling of a
hexagon with side lengths $2n$,
$2n$, $2n+1$, $2n$, $2n$, $2n+1$ contains
the (horizontal) rhombus with bottom-most vertex $(2n,2n)$ can be
written in the form
\begin{multline} \label{e3.1}
\text {SUM}(n):=\sum _{k=0} ^{n-1}
\frac {2n\,(2n+1)!} {(2 n+1)_{4 n}}
\binom {2 n} n\binom {3 n}n\,
 2^{2 n-2} \,(n+3/2)_k \,(2 n+1)_k \\
\cdot(n+k+2)_{n-k-1}
 \,(2 n+k+2)_{n-k-1}\frac {(1/2)_{n-k-1}} {(n-k-1)!}.
\end{multline}
Next we take it as an input for Zeilberger's algorithm (we are using
Zeilberger's {\sl Maple} implementation here):

\MATH
\leavevmode    |\MATHbackslash \MATHhoch /|     Maple V Release 4 (Uni Wien)
.\MATHtief |\MATHbackslash |   |/|\MATHtief . Copyright (c) 1981-1996 by Waterloo Maple Inc. All rights
\leavevmode \MATHbackslash   MAPLE  /  reserved. Maple and Maple V are registered trademarks of
\leavevmode \MATHkleiner \MATHtief \MATHtief \MATHtief \MATHtief  \MATHtief \MATHtief \MATHtief \MATHtief \MATHgroesser   Waterloo Maple Inc.
\leavevmode      |       Type ? for help.
\MATHgroesser  read ekhad:

\MATHgroesser  ezra(zeillim);
\leavevmode                       zeillim(SUMMAND,k,n,N,alpha,beta)
\leavevmode         Similar to zeil(SUMMAND,k,n,N) but outputs a recurrence for
\leavevmode                 the sum of SUMMAND from k=alpha to k=n-beta .
\leavevmode       Outputs the recurrence operator, certificate and right hand side.
\leavevmode       For example, "zeillim(binomial(n,k),k,n,N,0,1);" gives output of
\leavevmode                               N-2, k/(k-n-1),1
\leavevmode            which means that SUM(n):=2\MATHhoch n-1 satisfies the recurrence
\leavevmode               (N-2)SUM(n)=1, as certified by R(n,k):=k/(k-n-1)

\MATHgroesser  zeillim(2*n*(2*n+1)!/rf(2*n+1,4*n)*binomial(2*n,n)*binomial(3*n,n)*
\MATHgroesser  2\MATHhoch (2*n-2)*rf(n+3/2,k)*rf(2*n+1,k)*rf(n+k+2,n-k-1)*
\MATHgroesser  rf(2*n+k+2,n-k-1)*rf(1/2,n-k-1)/(n-k-1)!,k,n,N,0,1);

-1 + N, 1/6 (1 - 2 n + 2 k) (-288 n  - 432 n  k - 912 n  - 1440 n  k

\leavevmode            4  2        3  2         3           3        2  2        2
\leavevmode     - 216 n  k  - 414 n  k  - 1746 n  k - 1096 n  - 189 n  k  - 612 n

\leavevmode           2  3        2           3                           2
\leavevmode     + 36 n  k  - 907 n  k + 48 n k  - 152 n - 163 n k + 43 n k  - 12 + 4 k

\leavevmode           2       3    /
\leavevmode     + 32 k  + 16 k )  /  ((-n + k) n (6 n + 1) (6 n + 5) (2 n + k + 2)
\leavevmode                      /

\leavevmode                  2
\leavevmode    (1 + 3 n + 2 n )),

\leavevmode                  2               3   n      3       2
\leavevmode        GAMMA(3 n)  GAMMA(n + 1/2)  64  (36 n  + 60 n  + 29 n + 3)
\leavevmode    1/2 ----------------------------------------------------------
\leavevmode               2   3/2         3                      2
\leavevmode        (n + 1)  Pi    GAMMA(n)  (6 n + 5) (6 n + 1) n  GAMMA(6 n)

\MATHgroesser
\goodbreakpoint%
\endMATH
\noindent
It tells us that the expression $\text {SUM}(n)$ in \eqref{e3.1}
satisfies the recurrence
\begin{equation} \label{e3.2}
\text {SUM}(n+1)-\text {SUM}(n)=
\frac{\left( 3 + 29 n + 60 n^2 + 36 n^3 \right) \,
    {({ \textstyle 3 n-1})! }^2\,
    {({ \textstyle 2 n})! }^3}{2 n^2\,
    {\left( n+1 \right) }^2\,\left( 6 n +1\right) \,
    \left(  6 n+5 \right) \,{({ \textstyle n-1})! }^3\,
    ({ \textstyle 6 n-1})! \,
    {{ \textstyle n}! }^3}.
\end{equation}
(The first term in the output, {\tt -1 + N}, encodes the form of
the left-hand side of
\eqref{e3.2}, the third term gives the right-hand side. The middle
term is the so-called {\it certificate} which provides a proof of the
recurrence.)
So it just remains to check that the expression \eqref{e1.1} with
$f_{0,0,1,0,0}(n)$ as in \eqref{L1.6} satisfies the same recurrence
and agrees with \eqref{e3.1} for $n=1$,
which is of course a routine task.

On the other hand, the expression for $f_{2,1,1,2,1}(n)$ given in
\eqref{L1.2} cannot be established in the same way by appealing to
a special case of
one of (A)--(H). Still, it is also a theorem, thanks to the following
simple observation: suppose that we consider a hexagon with side
lengths $a$, $b$, $c$, $a$, $b$, $c$, where $a=b$, and a rhombus on
the horizontal symmetry axis of the hexagon. Let us imagine that this
rhombus were the one in Figure~3 (consisting of the triangles
labelled $0$ and $1$). Let us denote the probability that a randomly
chosen tiling contains this rhombus by $p$.
Since the rhombus is on the symmetry axis, the
probability that a randomly chosen tiling contains the rhombus
consisting of the triangles labelled $0$ and $2$ is equal to the
probability that it contains the rhombus consisting of the
triangles labelled $0$ and $3$. Let us denote this probability by
$q$. Any tiling must contain exactly one
of these three rhombi, hence we have $p+q+q=p+2q=1$. Therefore,
whenever there is a
single sum formula available for $p$, there is also one for
$q$. To come back to our example, the rhombus whose bottom-most point
has coordinates
$(2n+2,2n+1)$ in a hexagon with side lengths $2n+2$, $2n+1$, $2n+1$,
$2n+2$, $2n+1$, $2n+1$, can be seen as such a rhombus consisting of
triangles labelled $0$ and $3$, where the bottom-most point of the
rhombus consisting of the
triangles labelled $0$ and $1$ has coordinates $(2n+1,2n+1)$. This puts
us in Case~(E), with $a'=x'=-1$, (to see this one has to
replace $n$ by $n-1$ in the above coordinatization),
and thus the claimed expression for
$f_{2,1,1,2,1}(n)$ can be proved in the same manner as we proved
the expression for $f_{0,0,1,0,0}(n)$ above.

Again, for the sake of
completeness, we list the additional vectors $(a,b,c,x,y)$
for which single sums are available for $f_{a,b,c,x,y}(n)$ by the
above observation.

\begin{enumerate}
\item[(C')] By \cite[Theorem~1]{FuKrAD}:
$(2a',0,0,a'+1,a')$,
$(2a',0,0,a',a')$,
$(2a',0,0,a',a'-1)$,
$(2a',0,\break 0,a'-1,a'-1)$,
$(0,2a',0,a'+1,0)$,
$(0,2a',0,a',-1)$,
$(0,2a',0,a',0)$,
$(0,2a',0,a'-1,-1)$,
for an integer $a'$.
\item[(D')] By \cite[Theorem~2]{FuKrAD}:
$(2a'+1,1,1,a'+2,a'+1)$,
$(2a'+1,1,1,a'+1,a'+1)$,
$(2a'+1,1,1,a'+1,a')$,
$(2a'+1,1,1,a',a')$,
$(1,2a'+1,1,a'+2,1)$,
$(1,2a'+1,1,a'+1,0)$,
$(1,2a'+1,1,a'+1,1)$,
$(1,2a'+1,1,a',0)$,
for an integer $a'$.
\item[(E')] By \cite[Theorems~1 and 2]{FuKrAC}:
$(2a',0,0,a'+x',a'-x')$,
$(2a',0,0,a'+x',a'-x'-1)$,
$(2a'+1,0,0,a'+x',a'-x')$,
$(2a'+1,0,0,a'+x',a'-x'+1)$,
$(0,2a',0,a'+x',2x'-1)$,
$(0,2a',0,a'+x',2x')$,
$(0,2a'+1,0,a'+x',2x'-2)$,
$(0,2a'+1,0,a'+x',2x'-1)$,
for integers $a'$ and $x'$.
\item[(F')] By \cite[Theorems~1 and 2]{FuKrAC}:
$(2a',1,1,a'+x',a'-x'+1)$,
$(2a',1,1,a'+x',a'-x')$,
$(2a'+1,1,1,a'+x',a'-x'+2)$,
$(2a'+1,1,1,a'+x',a'-x'+1)$,
$(1,2a',1,a'+x',2x'-1)$,
$(1,2a',1,a'+x',2x')$,
$(1,2a'+1,1,a'+x',2x'-2)$,
$(1,2a'+1,1,a'+x',2x'-1)$,
for integers $a'$ and $x'$.
\item[(G')] By \cite[Theorem~3]{FuKrAD}:
$(2a'+1,0,0,a'+2,a'+2)$,
$(2a'+1,0,0,a'+1,a'+2)$,
$(2a'+1,0,0,a',a')$,
$(2a'+1,0,0,a'-1,a')$,
$(0,2a'+1,0,a'+2,0)$,
$(0,2a'+1,0,a'+1,-1)$,
$(0,2a'+1,0,a',0)$,
$(0,2a'+1,0,a'-1,-1)$,
for an integer $a'$.
\item[(H')] By \cite[Theorem~4]{FuKrAD}:
$(2a',1,1,a'+2,a'+2)$,
$(2a',1,1,a'+1,a'+2)$,
$(2a',1,1,a',a')$,
$(2a',1,1,a'-1,a')$,
$(1,2a',1,a'+2,1)$,
$(1,2a',1,a'+1,0)$,
$(1,2a',1,a',1)$,
$(1,2a',1,\break a'-1,0)$,
for an integer $a'$.
\item[(I')] By \cite[Theorem~5]{FuKrAD}:
$(2a',0,0,a'+2,a'+1)$,
$(2a',0,0,a'+1,a'+1)$,
$(2a',0,0,a'-1,a')$,
$(2a',0,0,a'-2,a')$,
$(0,2a',0,a'+2,0)$,
$(0,2a',0,a'+1,-1)$,
$(0,2a',0,a'-1,0)$,
$(0,2a',0,a'-2,-1)$,
for an integer $a'$.
\item[(J')] By \cite[Theorem~6]{FuKrAD}:
$(2a'+1,1,1,a'+3,a'+2)$,
$(2a'+1,1,1,a'+2,a'+2)$,
$(2a'+1,1,1,a',a'-1)$,
$(2a'+1,1,1,a'-1,a'-1)$,
$(1,2a'+1,1,a'+3,1)$,
$(1,2a'+1,1,a'+2,0)$,
$(1,2a'+1,1,a',1)$,
$(1,2a'+1,1,a'-1,0)$,
for an integer $a'$.
\end{enumerate}

Table~1 lists the special cases that we considered
in \eqref{L1.1}--\eqref{L1.3}, \eqref{L1.4}, \eqref{L1.5},
\eqref{L1.6a}--\eqref{L1.7}, for which proofs are available, together
with an indication from which of the Cases~(A)--(H), respectively
(C')--(H'), these proofs come from. In particular, the expression for
$f_{0,0,1,0,0}(n)$ in \eqref{L1.6} (together with \eqref{e1.1}) provides
the formula for the values of Pauling's bond order in Tabelle~2 of
\cite{JohPAA} for higher naphtalenes $N(p)$ of odd order $p=2n-1$, and
the expression for
$f_{1,1,0,1,0}(n)$ in \eqref{L1.5} (together with \eqref{e1.1}) provides
the formula for the values of Pauling's bond order in Tabelle~3 of
\cite{JohPAA} for higher pyrenes $P(p)$ of odd order $p=2n-1$. The values
in Tabelle~1 of \cite{JohPAA} for higher benzenes $B(p)$ are expressed
by \eqref{e1.1} with $a=b=c=x=y=0$ and $a=b=c=x=y=1$, respectively,
with the expression for $f_{0,0,0,0,0}(n)$ given in \eqref{L1.6a} and the one
for $f_{1,1,1,1,1}(n)$ given in \eqref{L1.3}. (The latter formulas have
already been stated in \cite[Cor.~7, (1.8) and (1.9)]{FuKrAD}. In
fact, Corollary~7 of \cite{FuKrAD} contains some more evaluations of
this kind.)

\begin{table}
\small
$$
\begin{array}{rrrrr|rl}
&&&&&\multicolumn{2}{c}{\text {entry in (A)--(H), } }\\
a,&b,&c,&x,&y& \multicolumn{2}{c}{\text{respectively (C')--(H')}} \\ \hline
-1,& -1,& 0,& -1,& -1& \text {(A), }& a'=b'=-1 \\ \hline
 2,& 2,& 1,& 2,& 1& \text{(B), }&  a'=b'=1 \\ \hline
2,& 1,& 1,& 2,& 1& \text{(F'), }& a'=x'=1 \\ \hline
 2,& 1,& 1,& 1,& 1& \text{(F'), }& a'=1,\ x'=0  \\ \hline
  1,& 2,& 1,& 2,& 1& \text{(F'), }& a'=x'=1  \\ \hline
 1,& 2,& 1,& 1,& 0& \text{(F'), }& a'=1,\ x'=0 \\ \hline
  -1,& 0,& 0,& 0,& -1& \text{(E'), }& a'=-1,\ x'=1  \\ \hline
 -1,& 0,& 0,& -1,& -1& \text{(E'), }& a'=-1,\ x'=0  \\ \hline
  0,& -1,& 0,& 0,& 0& \text{(E'), }& a'=-1,\ x'=1 \\ \hline
 0,& -1,& 0,& -1,& -1& \text{(E'), }& a'=-1,\ x'=0 \\ \hline
1,& 1,& 1,& 1,& 1& \text{(C), }& a'=b'=0 \\ \hline
1,& 1,& 0,& 1,& -1& \text{(H), }& a'=b'=0 \\ \hline
1,& 1,& 0,& 1,& 1&\text{(H), }& a'=b'=0 \\ \hline
1,& 1,& 1,& 0,& -1& \text{(F), }& a'=1,\ x'=0 \\ \hline
1,& 1,& 1,& 0,& 1& \text{(J'), }& a'=0 \\ \hline
1,& 1,& 1,& 1,& 0& \text{(C), }& a'=b'=0 \\ \hline
  0,& 2,& 0,& 1,& 0& \text{(E'), }& a'=1,\ x'=0 \\ \hline
2,& 0,& 0,& 0,& 1& \text{(I'), }& a'=1 \\ \hline
2,& 0,& 0,& 1,& 0& \text{(E'), }& a'=1,\ x'=0 \\ \hline
2,& 0,& 0,& 1,& 1& \text{(E'), }& a'=1,\ x'=0 \\ \hline
2,& 0,& 0,& 2,& 0& \text{(E'), }& a'=x'=1 \\ \hline
0,& 2,& 1,& 1,& 0& \text{(B), }& a'=0,\ b'=1 \\ \hline
1,& 1,& 0,& 1,& 0& \text{(A), }& a'=b'=0 \\ \hline
0,& 0,& 0,& 0,& 0& \text{(D), }& a'=b'=0 \\ \hline
0,& 0,& 1,& 0,& 0& \text{(B), }& a'=b'=0 \\ \hline
3,& 3,& 0,& 3,& 1& \text{(A), }& a'=b'=1
\end{array}
$$
\caption{}
\end{table}

\medskip
Coming back to the original goal, a proof of the Conjecture for
arbitrary $a$, $b$, $c$, $x$ and $y$, it may seem that it should be
at least possible to achieve this in the Cases~(A)--(H) and
(C')--(H'), where single sum formulas are available. For, for each
specific choice of $a$, $b$, $c$, $x$ and $y$ out of one of these
cases, an identity of the form ``single sum = closed form" has to be
proved. So one would try to follow the strategy that was suggested in
complete generality at the beginning of this section: apply some
manipulations (using hypergeometric transformation and summation
formulas, for example) until the desired expression is obtained.
This task is much less daunting here, since we are dealing now with a
single sum, not with a triple sum. Moreover, as it turns out,
the sums that occur are
very familiar objects in hypergeometric theory (we refer the reader to
\cite{BailAA,SlatAC,GaRaAA} for information on this theory),
they turn out to be
balanced $_4F_3$-series, respectively very-well-poised $_7F_6$-series.
(For example, the series in \eqref{e3.1} is a balanced $_4F_3$-series.)
For these series there are a lot of summation and transformation
formulas known. However, and this is somehow mysterious, I was
not able to establish any of the theorems that I presented here in this
classical manner (i.e., without the use of Zeilberger's algorithm),
not to mention a general theorem for an infinite family of
parameters. As already said at the beginning of this section,
the biggest stumbling block in such an attempt is the
question of how one would be able to isolate ``$1/3$" from the ``rest."
So, potentially, there is a hierarchy of interesting hypergeometric
identities lurking behind the scene which has not yet been
discovered.
\end{section}


\begin{thebibliography}{10}

\bibitem{BailAA}
W. N. Bailey, {\em Generalized hypergeometric series},
Cambridge University Press, Cambridge, 1935.

\bibitem{CiKrAA}
M.~Ciucu and C.~Krattenthaler, {\em The number of centered lozenge tilings of a
symmetric hexagon\/}, J. Combin\@. Theory Ser.~A {\bf 86} (1999), 103--126.

\bibitem{CoLPAA}
H.    Cohn, M. Larsen and J. Propp,
{\em The shape of a typical boxed plane partition},
New York J. Math\@. {\bf 4} (1998), 137--166.

\bibitem{DaToAA}
G. David and C. Tomei, {\em The problem of the calissons},
Amer\@. Math\@. Monthly\@. {\bf 96} (1989), 429--431.

\bibitem{FiscAA}
I. Fischer, {\em Enumeration of rhombus tilings of a hexagon which
contain a fixed rhombus in the centre},
J. Combin\@. Theory Ser.~A (to appear), {\tt math/9906102}.

\bibitem{FuKrAC}
M. Fulmek and C. Krattenthaler, {\em The number of rhombus
tilings of a symmetric hexagon which contain a fixed rhombus on the
symmetry axis, I}, Ann.\ Combin.\ 2 (1998), 19--40.

\bibitem{FuKrAD}
M. Fulmek and C. Krattenthaler, {\em The number of rhombus
tilings of a symmetric hexagon which contain a fixed rhombus on the
symmetry axis, II}, Europ\@. J. Combin\@. {\bf 21} (2000), 601--640.

\bibitem{GaRaAA} G.    Gasper and M. Rahman,
{\em Basic Hypergeometric Series},
Encyclopedia of Mathematics And Its Applications~35, Cambridge
  University Press, Cambridge, 1990.

\bibitem{HeGeAA}
H.    Helfgott and I. M. Gessel,
{\em Exact enumeration of tilings of diamonds and hexagons with
defects}, Electron\@. J. Combin\@. {\bf 6} (1) (1999), \#R16, 26~pp.

\bibitem{JohaAB}
K.    Johansson, {\em Nonintersecting paths,
random tilings and random matrices}, preprint,
{\tt math/0011250}.

\bibitem{JohPAA}
P.    John, {\em
\"Uber ein einfaches Wachstum hexagonaler Systeme und das Verhalten der
Paulingschen Bindungsordnung},
Wiss\@. Zeitschr\@. FSU Jena {\bf 39} (1990),
192--200.

\bibitem{KupeAA}
G.    Kuperberg, {\em Symmetries of plane partitions and the permanent
  determinant method}, J.~Combin\@. Theory Ser\@.~A {\bf 68} (1994),
115--151.

\bibitem{KratBN}
C.    Krattenthaler,
{\em Advanced determinant calculus},
S\'eminaire Lotharingien Combin\@. {\bf 42} \rm(``The Andrews Festschrift")
(1999), paper~B42q, 67~pp.

\bibitem{MacMAA}
P.~A.~MacMahon, {\em Combinatory Analysis\/}, vol.~2, Cambridge University
Press, 1916; reprinted by Chelsea, New York, 1960.

\bibitem{PaScAA} P.    Paule and M. Schorn, {\em A
Mathematica version of Zeilberger's algorithm for proving binomial
coefficient identities}, J. Symbol\@. Comp\@. {\bf 20} (1995),
673--698.

\bibitem{PeWZAA} M.    Petkov\v sek, H. Wilf and D. Zeilberger,
{\em A=B}, A.~K.~Peters, Wellesley, 1996.

\bibitem{PropAA}
J.~Propp, {\em Twenty open problems on enumeration of matchings\/},
manuscript, 1996, {\tt math/9801060}.

\bibitem{PropAH}
J.~Propp, {\em Enumeration of matchings: Problems and progress},
in: "New Perspectives in Algebraic
 Combinatorics", L.~Billera, A.~Bj\"orner, C.~Greene, R.~Simion, and
R.~P.~Stanley, eds.,
 Mathematical Sciences Research Institute Publications, vol.~38,
Cambridge University Press, 1999, pp.~255--291.

\bibitem{SlatAC}
L.~J.~Slater,
{\em Generalized hypergeometric functions\/},
Cambridge University Press, Cambridge, 1966.

\bibitem{VerbAA}
P.    Verbaeten, {\em Rekursiebetrekkingen voor
lineaire hypergeometrische funkties},
Proefschrift voor het doctoraat in
de toegepaste wetenschapen, Katholieke Universiteit te Leuven,
Heverlee, Belgium, 1976.

\bibitem{WegsAA}
K.    Wegschaider, {\em Computer generated
proofs of binomial multi-sum identities}, diploma thesis,
Johannes Kepler University, Linz, Austria, 1997; available from
{\tt http://www.risc.uni-linz.ac.at/research/combinat/risc/publications}.

\bibitem{WiZeAC}
H. S. Wilf and D. Zeilberger, {\em An
algorithmic proof theory for hypergeometric (ordinary and ``$q$")
multisum/integral identities}, Invent\@. Math\@. {\bf 108} (1992),
575--633.

\bibitem{ZeilAM} D.    Zeilberger,
{\em A fast algorithm for proving terminating hypergeometric
identities}, Discrete Math\@. {\bf 80} (1990), 207--211.

\bibitem{ZeilAV} D.    Zeilberger,
{\em The method of creative telescoping},
J. Symbolic Comput\@. {\bf 11} (1991), 195--204.

\end{thebibliography}
\end{document}